\documentclass[12pt,a4paper]{article}
\usepackage{amsfonts,mathtools,etoolbox,amsthm,amssymb,amsmath,bbm,enumitem,nicefrac,xcolor,xparse,hyperref}

\usepackage[utf8]{inputenc}
\usepackage[english]{babel}
\usepackage{float} % Place figures exactly where you want with \begin{figure}[H]
\usepackage[T1]{fontenc} %font stuff

\usepackage[nocompress]{cite}
\usepackage[capitalise,sort]{cleveref}
\newcommand{\B}{\mathbb{B}}

\crefname{enumi}{item}{items}
\crefname{equation}{}{}
\numberwithin{equation}{section}
\crefname{figure}{Figure}{Figures}
\crefname{listing}{Source code}{Source codes}
\crefname{lstlisting}{Source code}{Source codes}
\crefname{cor}{Corollary}{Corollaries}
\crefname{subsection}{Subsection}{Subsections}

\hypersetup{
	colorlinks,
	linkcolor={blue!80!black},
	citecolor={green},
	urlcolor={blue!80!black}
}

\usepackage[margin=1in]{geometry}

\theoremstyle{plain}
\newtheorem{theorem}{Theorem}[section]

\newtheorem{lemma}[theorem]{Lemma}
\newtheorem{prop}[theorem]{Proposition}
\newtheorem{cor}[theorem]{Corollary}
\theoremstyle{definition}
\newtheorem{definition}[theorem]{Definition}
\newtheorem{remark}[theorem]{Remark}

\newcommand{\N}{\mathbb{N}}
\newcommand{\R}{\mathbb{R}}
\newcommand{\C}{\mathbb{C}}

\newcommand{\eps}{\varepsilon}
\newcommand{\fl}{\mathcal{L}}
\newcommand{\cK}{\mathcal{K}}

\newcommand{\cL}{\mathcal{L}}

\newcommand{\qqandqq}{\qquad \text{and} \qquad}

\newcommand{\sr}{\rho}
\newcommand{\fd}{d}

\newcommand{\cond }{ \chi }

\newcommand{\m}{\mathbf{m}}
\newcommand{\bbM}{\mathbb{M}}
\newcommand{\lp}{\psi}
\newcommand{\mlr}{\overline{\gamma}}
\newcommand{\vone}{\mathbf{1}}
\newcommand{\dia}[1]{\operatorname{diag}( #1 )}
\newcommand{\unit}[1]{\operatorname{I}_{#1}}
\newcommand{\spec}[1]{\Lambda \rbr*{ #1 } }
\newcommand{\Hess}{\operatorname{Hess}}
\newcommand{\sprod}{\textstyle\prod}
\newcommand{\cp}{\odot}

\newcommand{\const}{\mathfrak{c}}
\newcommand{\with}{\curvearrowleft}

\newcommand{\indicator}[1]{\mathbbm{1}_{\smash{#1}}}

\DeclarePairedDelimiter{\norm}{\lVert}{\rVert}
\DeclarePairedDelimiter{\abs}{\lvert}{\rvert}
\DeclarePairedDelimiter{\rbr}{(}{)}
\DeclarePairedDelimiter{\br}{[}{]}
\DeclarePairedDelimiter{\cu}{\{}{\}}
\DeclarePairedDelimiter{\spro}{\langle}{\rangle}

\ExplSyntaxOn

\bool_new:N \g_noteobserve

\NewDocumentCommand{\setnote}{}{
	\bool_gset_true:N \g_noteobserve
}

\NewDocumentCommand{\setobserve}{}{
	\bool_gset_false:N \g_noteobserve
}

\NewDocumentCommand{\nobs}{ o }{
	\IfValueT{#1}{
		\str_if_eq:noTF {note} {#1} {
			\bool_gset_true:N \g_noteobserve
		} {
			\str_if_eq:noTF {Note} {#1} {
				\bool_gset_true:N \g_noteobserve
			} {
				\bool_gset_false:N \g_noteobserve
			}
		}
	}
	\bool_if:nTF { \g_noteobserve } {
		\bool_gset_false:N \g_noteobserve
		note
	} {
		\bool_gset_true:N \g_noteobserve
		observe
	}
	\IfValueF{#1}{~}
}

\NewDocumentCommand{\Nobs}{ o }{
	\IfValueT{#1}{
		\str_if_eq:noTF {note} {#1} {
			\bool_gset_true:N \g_noteobserve
		} {
			\str_if_eq:noTF {Note} {#1} {
				\bool_gset_true:N \g_noteobserve
			} {
				\bool_gset_false:N \g_noteobserve
			}
		}
	}
	\bool_if:nTF { \g_noteobserve } {
		\bool_gset_false:N \g_noteobserve
		Note
	} {
		\bool_gset_true:N \g_noteobserve
		Observe
	}
	\IfValueF{#1}{~}
}

\seq_new:N \g_cflist_loaded
\seq_new:N \g_cflist_pending

\NewDocumentCommand{\cfadd}{ m }
{
	\seq_if_in:NnF \g_cflist_loaded { #1 } {
		\seq_if_in:NnF \g_cflist_pending { #1 } {
			\seq_gput_right:Nn \g_cflist_pending { #1 }
		}
	}
}

\NewDocumentCommand{\cfconsiderloaded}{ m }{
	\seq_gput_right:Nn \g_cflist_loaded {#1}
}

\NewDocumentCommand{\cfremove}{ m }
{
	\seq_gremove_all:Nn \g_cflist_pending { #1 }
}

\NewDocumentCommand{\cfload}{ o }
{
	\seq_if_empty:NTF \g_cflist_pending {\unskip} {
		(cf.\ \cref{\seq_use:Nn \g_cflist_pending {,}})\IfValueTF{#1}{#1~}{\unskip}
		\seq_gconcat:NNN \g_cflist_loaded \g_cflist_loaded \g_cflist_pending
		\seq_gclear:N \g_cflist_pending
	}
}

\NewDocumentCommand{\cfclear} {} {
	\seq_gclear:N \g_cflist_loaded
	\seq_gclear:N \g_cflist_pending
}

\NewDocumentCommand{\cfout}{ o }
{
	\seq_if_empty:NTF \g_cflist_pending {\unskip} {
		(cf.\ \cref{\seq_use:Nn \g_cflist_pending {,}})\IfValueTF{#1}{#1~}{\unskip}
		\seq_gclear:N \g_cflist_pending
	}
}

\NewDocumentCommand{\ifnocf} { m } {
	\seq_if_empty:NT \g_cflist_pending { #1 }
}

\ExplSyntaxOff

\NewDocumentEnvironment{cproof}{m}
{\begin{proof}[Proof of \cref{#1}]}%
	{\noindent The proof of \cref{#1} is thus complete.
\end{proof}}

\NewDocumentEnvironment{cproof2}{m}
{\begin{proof}[Proof of \cref{#1}]}%
	{\noindent This completes the proof of \cref{#1}.
\end{proof}}

\ExplSyntaxOn

\bool_new:N \g_forexample

\NewDocumentCommand{\eg}{ o }{
	\IfValueT{#1}{
		\str_if_eq:noTF {fe} {#1} {
			\bool_gset_true:N \g_forexample
		} {\bool_gset_false:N \g_forexample}
	}
	\bool_if:nTF { \g_forexample } {
		\bool_gset_false:N \g_forexample
		for~example
	}{
		\bool_gset_true:N \g_forexample
		for~instance
	}
}

\ExplSyntaxOff

\makeatletter
\ExplSyntaxOn
\seq_new:N \g_abbrs
\prop_new:N \g_abbr_counts
\tl_new:N \l_abbr_count_tl

\NewDocumentCommand{\abbr}{m m O{#1} m m O{#4}}{
	\expandafter\newcommand\csname#3\endcsname[1][]{
		\seq_if_in:NnTF \g_abbrs {#1} {
			\prop_get:NnN \g_abbr_counts {#1} \l_abbr_count_tl
			\prop_gput:Nnx \g_abbr_counts {#1} {\int_eval:n {\l_abbr_count_tl + 1}}
			\hyperref[#1]{#1}
		} {
			\seq_gput_left:Nn \g_abbrs {#1}
			\prop_gput:Nnn \g_abbr_counts {#1} {1}
			\expandafter\gdef\csname#1@def\endcsname{#2}
			\phantomsection\label{#1}
			\str_if_eq:nnTF{##1}{}{\emph{#2}}{##1}~(\hyperref[#1]{#1})
		}
	}
	\expandafter\newcommand\csname#6\endcsname[1][]{
		\seq_if_in:NnTF \g_abbrs {#1} {
			\prop_get:NnN \g_abbr_counts {#1} \l_abbr_count_tl
			\prop_gput:Nnx \g_abbr_counts {#1} {\int_eval:n {\l_abbr_count_tl + 1}}
			\hyperref[#1]{#4}
		} {
			\expandafter\gdef\csname#1@def\endcsname{#5}
			\seq_gput_left:Nn \g_abbrs {#1}
			\prop_gput:Nnn \g_abbr_counts {#1} {1}
			\phantomsection\label{#1}
			\str_if_eq:nnTF{##1}{}{\emph{#5}}{##1}~(\hyperref[#1]{#4})
		}
	}
}

\ExplSyntaxOff
\makeatother

\abbr{ReLU}{rectified linear unit}{ReLUs}{rectified linear units}
\abbr{OP}{optimization problem}{OPs}{optimization problems}
\abbr{SOP}{stochastic optimization problem}{SOPs}{stochastic optimization problems}
\abbr{ANN}{artificial neural network}{ANNs}{artificial neural networks}
\abbr{DL}{deep learning}{DLs}{deep learnings}
\abbr{iid}{independent and identically distributed}{iids}{deep learnings}
\abbr{GD}{gradient descent}{GDs}{gradient descent}
\abbr{SGD}{stochastic gradient descent}{SGDs}{stochastic gradient descent}
\abbr{Adam}{adaptive moment estimation}{Adams}{stochastic gradient descent}
\abbr{RMSprop}{root mean square propagation}{RMSprops}{root mean square propagations}
\abbr{AdamW}{\Adam\ with weight decay}{AdamWs}{\Adams with weight decay}
\abbr{EMA}{exponential moving average}{EMAs}{exponential moving averages}
\abbr{RPA}{Ruppert--Polyak averaging}{RPAs}{Ruppert--Polyak averagings}
\abbr{ODE}{ordinary differential equation}{ODEs}{ordinary differential equations}
\abbr{Adagrad}{adaptive gradient}{Adagrads}{adaptive gradients}

\title{Sharp higher order convergence\\ rates for the Adam optimizer}

\author{Steffen Dereich$^{1}$, Arnulf Jentzen$^{2,3}$, and Adrian Riekert$^{4}$\bigskip\\
	\small{$^1$ Institute for Mathematical Stochastics, University of M\"unster,}\vspace{-0.1cm}\\
	\small{Germany; e-mail: \texttt{steffen.dereich}\textcircled{\texttt{a}}\texttt{uni-muenster.de}}\smallskip\\
	\small{$^2$ School of Data Science and Shenzhen Research Institute of Big Data,}\vspace{-0.1cm}\\
	\small{The Chinese University of Hong Kong, Shenzhen (CUHK-Shenzhen),}\vspace{-0.1cm}\\
	\small{China; e-mail: \texttt{ajentzen}\textcircled{\texttt{a}}\texttt{cuhk.edu.cn}}\smallskip\\
	\small{$^3$ Applied Mathematics: Institute for Analysis and Numerics,}\vspace{-0.1cm}\\
	\small{University of M\"unster, Germany; e-mail: \texttt{ajentzen}\textcircled{\texttt{a}}\texttt{uni-muenster.de}}\smallskip\\
	\small{$^4$ Applied Mathematics: Institute for Analysis and Numerics,}\vspace{-0.1cm}\\
	\small{University of M\"unster, Germany; e-mail: \texttt{ariekert}\textcircled{\texttt{a}}\texttt{uni-muenster.de}}}

\date{\today}

\begin{document}
	
	\maketitle

	\begin{abstract}
		Gradient descent based optimization methods
		are the methods of choice to train deep neural networks
		in machine learning.
		Beyond the standard gradient descent method,
		also suitable modified variants of
		standard gradient descent involving
		acceleration techniques such as the momentum method
		and/or adaptivity techniques such as the RMSprop method
		are frequently considered optimization methods.
		These days the most popular of such sophisticated
		optimization schemes is presumably
		the \emph{Adam optimizer} that has been proposed in 2014 by Kingma and Ba.
		A highly relevant topic of research is to investigate
		the speed of convergence of such optimization methods.
		In particular, in 1964 Polyak showed that
		the standard gradient descent method
		converges in a neighborhood of a strict local minimizer
		with rate
		$ ( \cond - 1 ) ( \cond + 1 )^{ - 1 } $
		while momentum achieves the (optimal)
		strictly faster convergence rate
		$ ( \sqrt{ \cond } - 1 ) ( \sqrt{ \cond } + 1 )^{ - 1 } $
		where $ \cond \in (1,\infty) $
		is the condition number
		(the ratio of the largest and the smallest eigenvalue)
		of the Hessian of the objective function at the local minimizer.
		It is the key contribution of this work to
		reveal that Adam
		also converges with the \emph{strictly faster convergence rate}
		$ ( \sqrt{ \cond } - 1 ) ( \sqrt{ \cond } + 1 )^{ - 1 } $
		while RMSprop only converges with the convergence rate
		$ ( \cond - 1 ) ( \cond + 1 )^{ - 1 } $.
	\end{abstract}

	\tableofcontents
	
	\section{Introduction}

	In deep learning, one typically trains a class of
	deep \ANNs\ by means of a \GD\ based optimization method.
	Beside the standard \GD\ method \cite{Cauchy1847},
	often suitable modified variants
	of the standard \GD\ method
	such as \GD-type optimization
	methods involving acceleration techniques
	as the momentum \GD\ optimizer
	\cite{Polyak1964}
	or \GD-type optimization methods
	involving adaptive learning rate
	schedules as in the \Adagrad\ \cite{DuchiHazanSinger11}
	and the \RMSprop\ \cite{HintonSlides} methods
	are employed to solve optimization problems.
	Nowadays the most widespread optimization
	algorithm in deep learning is presumably
	the famous \Adam\ optimizer proposed
	in 2014 by Kingma and Ba~\cite{KingmaBa2014_Adam} that consists,
	up to marginal modifications, of a combination of
	the momentum and
	the \RMSprop\ optimizers
	(cf., \eg, \cite{ruder2017overviewgradientdescentoptimization,JentzenBookDeepLearning2025}).

	Due to the enormous importance of such gradient based optimization methods,
	it is a very active and important field of research to investigate
	the convergence properties of such optimization algorithms (cf., \eg,
	\cite{MR2142598,Bach2024_book_MIT_press,JentzenBookDeepLearning2025} and the references therein).
	The optimal convergence rates for both the
	standard \GD\ optimization method as well as
	the momentum \GD\ optimization method
	have already been known for a long time.
	In particular, in 1964
	Polyak~\cite[Theorem 9]{Polyak1964}
	(cf.\ also \cite{Polyak1963} and, \eg,
	\cite[Item (i) in Proposition~6.3.26]{JentzenBookDeepLearning2025})
	proved in the situation
	where $ d \in \N = \{ 1, 2, 3, \dots \} $
	is the dimension of the optimization problem,
	where
	$ \cL \colon \R^d \to \R $
	is the objective function
	(the function that one intends to minimize)
	in the optimization problem
	\begin{equation}
		\label{eq:optimization_problem}
		\min_{ \theta \in \R^d }
		\cL( \theta )
		,
	\end{equation}
	where $ \lp \in \R^d $
	is a local minimizer of $ \cL $,
	and where
	$ \cond \in (1,\infty) $
	is the \emph{condition number}
	(the ratio of the largest and the smallest eigenvalue)
	of the Hessian
	$ ( \Hess \cL )( \lp ) $
	of the objective function
	$ \cL $ at the local minimizer $ \lp $
	that there exists a learning rate
	(a step size)
	such that
	the standard \GD\ method
	converges in a neighborhood of
	the local minimizer $ \lp $
	with rate
	$ ( \cond - 1 ) ( \cond + 1 )^{ - 1 } $
	to the local minimizer $ \lp $
	while momentum \GD\ achieves
	the strictly faster convergence rate
	$
	( \sqrt{ \cond } - 1 ) ( \sqrt{ \cond } + 1 )^{ - 1 }
	$ (see also \cref{theo:item:sgd,theo:item:mom-sgd}
	in \cref{theo:intro:rate} below).
	We refer to \cref{eq:convergence_rate_comparison}
	below for the obvious fact
	that the momentum \GD\ convergence rate is
	indeed strictly faster than the standard \GD\ convergence
	rate. We also refer to \cref{sec:literature} below
	for a more comprehensive
	% short
	literature overview on the mathematical analysis of
	such gradient based optimization methods.

	While the optimal rates of convergence
	of the standard \GD\ and
	the momentum \GD\ methods have been
	known for more than half a century,
	it remained an open problem of research
	to establish optimal convergence rates
	for the \Adam\ optimizer,
	presumably the most popular optimizer
	in deep learning.
	It is precisely the subject of this
	article to close this gap and to establish
	optimal convergence rates for the
	famous \Adam\ and the \RMSprop\ optimizers
	which are both based on
	\emph{adaptive learning rate schedules}.
	In \cref{theo:intro:rate}
	in \cref{sec:main_result}
	below we briefly
	outline the contribution of this article
	within this introductory section.
	In our formulation of \cref{theo:intro:rate}
	we employ the concept of a suitable
	convergence rate.
	This conceptuality is presented
	in \cref{def:convergence_order_strong}
	in \cref{sec:rate_of_convergence} below.

	\subsection{Rate of convergence}
	\label{sec:rate_of_convergence}

	In the following notion,
	\cref{def:convergence_order_strong} below,
	we aim to measure/describe the speed
	with which a certain general multi-step
	gradient based optimization method (see \cref{eq:optimization_method} below)
	converges
	to a local minimizer $ \lp \in \R^d $
	of the optimization problem
	in \cref{eq:optimization_problem}.

	\begin{samepage}
		\begin{definition}[Rate of convergence]
			\label{def:convergence_order_strong}
			Let $ \fd \in \N $, $ \varrho \in [0,1) $,
			$ \cond \in [1,\infty) $
			and let
			$ \Phi_n \colon ( \R^d )^{ n } \to \R^d $,
			$ n \in \N $,
			be functions.
			Then we say that
			$ ( \Phi_n )_{ n \in \N } $ converges
			with $ \cond $-rate $ \varrho $
			if and only if
			it holds\footnote{\Nobs that for all $ n \in \N $ it holds that
				$ \unit{n} \in \R^{ n \times n } $ is the identity matrix of size $ n $
				in the sense that for all $ v \in \R^n $ it holds that
				$
				\unit{n}
				v = v
				$ and \nobs that for all $ n \in \N $, $ A, B \in \R^{ n \times n } $ it holds that
				$ A \preceq B $ if and only if it holds for all $ v \in \R^n $ that
				$ \left< v, A v \right> \leq \left< v, B v \right> $.}
			for all
			$ \rho \in ( \varrho, 1 ) $,
			$
			\cL \in C^2( \R^d, \R )
			$,
			$ \lp \in \R^d $,
			$ r \in (0,\infty) $
			with
			$ (\nabla \cL ) ( \lp ) = 0 $
			and
			$
			r \unit{d} \preceq
			( \Hess \cL ) ( \lp )
			\preceq
			r \cond \unit{d}
			$
			that
			there exist $ \gamma, \delta \in (0,\infty) $ such that for all
			$ \Theta = ( \Theta_n )_{ n \in \N } \colon \N \to \R^{ \fd } $
			with
			$ \| \Theta_0 - \lp \| < \delta $
			and
			\begin{equation}
				\label{eq:optimization_method}
				\forall \, n \in \N \colon
				\Theta_n
				=
				\Theta_{ n - 1 }
				% \\
				-
				\gamma \Phi_n\bigl(
				% 		\Theta_0, \Theta_1, \dots, \Theta_{ n - 1 },
				( \nabla \cL )( \Theta_0 ) ,
				( \nabla \cL )( \Theta_1 ) ,
				\dots ,
				( \nabla \cL )( \Theta_{ n - 1 } )
				\bigr)
			\end{equation}
			it holds that
			$
			\sup_{ n \in \N }
			(
			\rho^{ - n } \| \Theta_n - \lp \|
			)
			< \infty
			$.
		\end{definition}
	\end{samepage}

	In \cref{def:convergence_order_strong}
	the natural number $ d \in \N $
	is the dimension of the optimization problem,
	the real number $ \varrho \in [0,1) $
	is the convergence rate,
	the parameter $ \cond \in [1,\infty) $
	is the condition number of the
	Hessian $ ( \Hess \cL )( \lp ) $
	of the objective function $ \cL \colon \R^d \to \R $
	at the local minimum point $ \lp $,
	and the functions
	$ \Phi_n \colon ( \R^d )^n \to \R^d $,
	$ n \in \N $,
	in \cref{eq:optimization_method}
	describe the gradient based
	optimization method under consideration.

	The standard \GD, the momentum \GD,
	the \RMSprop, the \Adam, as well as many
	other deep learning optimizers
	can be described through such a
	sequence of functions
	$ \Phi_n \colon ( \R^d )^n \to \R^d $,
	$ n \in \N $ (cf., \eg, \cite[Subsections~6.4 and 6.5]{MR4688424} and \cite[Section~6]{DoHannibalJentzen2024arXiv}).
	For instance, in the case of
	standard \GD\ we have for all $ n \in \N $,
	$ g_1, g_2, \dots, g_n \in \R^d $ that
	$
	\Phi_n( g_1, g_2, \dots, g_n ) = g_n
	$.
	We refer to
	\cref{def:GD}, \cref{def:RMSprop},
	\cref{def:momentum}, and \cref{def:Adam}
	in \cref{sec:gradient_descent,sec:momentum,sec:Adam} below
	for the formulation
	of the standard \GD, the momentum \GD, the \RMSprop,
	and the \Adam\ optimizers
	using the general functions
	$ ( \Phi_n )_{ n \in \N } $
	in \cref{eq:optimization_method}
	in \cref{def:convergence_order_strong} above.

	\subsection{Main result:
		Convergence speeds of GD, momentum,
		RMSprop, \& Adam}
	\label{sec:main_result}

	After having
	presented \cref{def:convergence_order_strong}
	above, we are now in the position
	to present in \cref{theo:intro:rate} below
	one of the main convergence
	results of this work in which we specify
	the convergence speeds
	of the standard \GD, the momentum \GD,
	the \RMSprop, and the \Adam\ optimizers.

	\cfclear
	\begin{samepage}
		\begin{theorem}[Convergence speeds of \GD, momentum, \RMSprop, and \Adam]
			\label{theo:intro:rate}
			\cfadd{def:convergence_order_strong,def:GD,def:momentum,def:RMSprop,def:Adam}
			Let $ \fd \in \N $,
			$ \cond \in [1,\infty) $,
			$ \eps \in (0,\infty) $, $ \beta \in (0,1) $.
			Then
			\begin{enumerate}[label=(\roman*)]
				
				\item
				\label{theo:item:sgd}
				for every \GD\ optimizer
				$ \Phi_n \colon ( \R^d )^{ n } \to \R^d $, $ n \in \N $,
				it holds that $ ( \Phi_n )_{ n \in \N } $
				converges with $ \cond $-order
				\textcolor{magenta}{$ \frac{ \cond - 1 }{ \cond + 1 } $},
				
				\item
				\label{theo:item:mom-sgd}
				there exists $ \alpha \in (0,1) $
				such that for every $ \alpha $-momentum \GD\ optimizer
				$ \Phi_n \colon ( \R^d )^{ n } \to \R^d $, $ n \in \N $,
				it holds that
				$ ( \Phi_n )_{ n \in \N } $
				converges with $ \cond $-order
				\textcolor{magenta}{$ \frac{ \sqrt{ \cond } - 1 }{ \sqrt{ \cond } + 1 } $},
				
				\item
				\label{theo:item:rmsprop}
				for every $ \beta $-$ \varepsilon $-\RMSprop\ optimizer
				$ \Phi_n \colon ( \R^d )^{ n } \to \R^d $, $ n \in \N $,
				it holds that
				$ ( \Phi_n )_{ n \in \N } $
				converges with $ \cond $-order
				\textcolor{magenta}{$ \frac{ \cond - 1 }{ \cond + 1 } $},
				and
				
				\item
				\label{theo:item:adam}
				there exists $ \alpha \in (0, 1 ) $
				such that for every $ \alpha $-$ \beta $-$ \varepsilon $-\Adam\ optimizer
				$ \Phi_n \colon ( \R^d )^{ n } \to \R^d $, $ n \in \N $,
				it holds that
				$ ( \Phi_n )_{ n \in \N } $
				converges with $ \cond $-order
				\textcolor{magenta}{$ \frac{ \sqrt{ \cond } - 1 }{ \sqrt{ \cond } + 1 } $}
			\end{enumerate}
			\cfload.
		\end{theorem}
	\end{samepage}

	\Cref{theo:item:sgd} in \cref{theo:intro:rate}
	is a consequence of \cref{cor:rmsprop_local}
	in \cref{ssec:convergence_RMSprop} below
	(applied for every
	$ \cond \in (1,\infty) $,
	$ r \in (0,\infty) $
	with
	$ \eps \with 1 $,
	$ \kappa \with r $,
	$ \cK \with r \cond $,
	$ \beta \with 0 $
	in the notation
	of \cref{cor:rmsprop_local}),
	\cref{theo:item:mom-sgd} in \cref{theo:intro:rate}
	is a consequence of \cref{cor:momentum-local-c2}
	in \cref{ssec:Adam_optimal_rate} below
	(applied for every
	$ \cond \in (1,\infty) $,
	$ r \in (0,\infty) $
	with
	$
	  \kappa \with r
	$,
	$
	  \cK \with r \cond
	$,
	$
  	  \alpha \with
	  ( \sqrt{ \cond } - 1 )
	  ( \sqrt{ \cond } + 1 )^{ - 1 }
	$,
	$
  	  \gamma \with r^{ - 1 / 2 } ( r \cond )^{ - 1 / 2 }
	$
	in the notation of \cref{cor:momentum-local-c2}),
	\cref{theo:item:rmsprop} in \cref{theo:intro:rate}
	is a consequence of \cref{cor:rmsprop_local}
	in \cref{ssec:convergence_RMSprop} below
	(applied for every
	$ \cond \in (1,\infty) $,
	$ r, \eps \in (0,\infty) $,
	$ \beta \in (0,1) $
	with
	$ \eps \with \eps $,
	$ \kappa \with r $,
	$ \cK \with r \cond $,
	$ \beta \with \beta $
	in the notation of \cref{cor:rmsprop_local}),
	and \cref{theo:item:adam} in \cref{theo:intro:rate}
	is a consequence of \cref{cor:adam-local-c2}
	in \cref{ssec:Adam_optimal_rate} below
	(applied for every
	$ \cond \in (1,\infty) $,
	$ r, \eps \in (0,\infty) $,
	$ \beta \in (0,1) $
	with
	$
	  \kappa \with r
	$,
	$
	  \cK \with r \cond
	$,
	$
	  \alpha \with
	  ( \sqrt{ \cond } - 1 )
	  ( \sqrt{ \cond } + 1 )^{ - 1 }
	$,
	$
	  \beta \with \beta
	$,
	$
 	  \gamma \with
 	  r^{ - 1 / 2 } ( r \cond )^{ - 1 / 2 } \eps
	$,
	$ \eps \with \eps $
	in the notation of \cref{cor:adam-local-c2}).

	It should, however, be pointed out
	that
	\cref{theo:item:sgd,theo:item:mom-sgd} in \cref{theo:intro:rate}
	have essentially been established by
	Polyak~\cite[Theorem 9]{Polyak1964}
	(cf.\ also Polyak~\cite{Polyak1963}
	and, \eg, \cite[Item (i) in Proposition~6.3.26]{JentzenBookDeepLearning2025})
	and are \emph{not} new
	contributions of this work.
	Instead the \emph{new contributions}
	of this work are
	to establish the convergence rates
	for the adaptive learning gradient based methods
	\RMSprop\ and \Adam, specifically,
	to reveal that \RMSprop\ converges
	with the standard order as the standard \GD\ method
	(see \cref{theo:item:rmsprop}
	in \cref{theo:intro:rate}) and
	to reveal
	that \Adam\ converges with
	the higher order
	as the accelerated momentum \GD\ method
	(see \cref{theo:item:adam} in \cref{theo:intro:rate}).
	We also refer to \cref{sec:literature} below
	for a more comprehensive literature overview
	on convergence analyses of gradient based
	optimization methods.

	We also note that
	for every condition number
	$ \cond \in (1,\infty) $
	we have that
	the convergence speed
	$
	(
	\sqrt{ \cond } - 1
	)
	(
	\sqrt{ \cond } + 1
	)^{ - 1 }
	$
	in \cref{theo:item:mom-sgd,theo:item:adam}
	in \cref{theo:intro:rate}
	is strictly faster than
	the convergence speed
	$
	( \cond - 1 )
	( \cond + 1 )^{ - 1 }
	$
	in \cref{theo:item:sgd,theo:item:rmsprop}
	in \cref{theo:intro:rate}.
	Indeed, we observe that
	for all $ \cond \in (1,\infty) $
	it holds that
	\begin{equation}
		\label{eq:convergence_rate_comparison}
		\frac{
			\sqrt{ \cond } - 1
		}{
			\sqrt{ \cond } + 1
		}
		=
		\frac{
			( \sqrt{ \cond } - 1 )
			( \sqrt{ \cond } + 1 )
		}{
			( \sqrt{ \cond } + 1 )^2
		}
		=
		\frac{
			\cond - 1
		}{
			( \sqrt{ \cond } + 1 )^2
		}
		=
		\frac{
			\cond - 1
		}{
			\kappa + 2 \sqrt{ \cond } + 1
		}
		<
		\frac{ \cond - 1 }{ \cond + 1 }
	\end{equation}
	(cf., \eg, \cite[Lemma~6.3.27 and Corollary~6.3.28]{JentzenBookDeepLearning2025}).
	We also note that, to some extend,
	the higher convergence order
	$
	  ( \sqrt{ \cond } - 1 ) ( \sqrt{ \cond } + 1 )^{ - 1 }
	$
	achieved by the momentum \GD\ and the \Adam\ optimizers
	is
	optimal and can
	not be improved within a large class optimization
	schemes; cf., \eg, \cite[Corollary~3 in Section~4]{MR4343730}
	and \cite[Theorem~2.1.13]{MR2142598}.

	\subsection{Literature review}
	\label{sec:literature}

	In this subsection we provide a short review
	on selected articles in the scientific literature
	that analyze gradient based optimization methods, particularly, adaptive gradient based optimization methods such as \Adam\ and \RMSprop.

    First, we note that for stochastic \OPs\ one
    can in general not expect a better/lower estimate
    for the
    \emph{strong optimization error (distance to the critical point)}
    than a constant
    multiplied by the square root
    of the learning rate in the sense that
    for every $ n \in \N $ we have that the
    root mean square distance between
    the optimization process at time $ n \in \N $
    and the minimizer of the stochastic \OP\ is bounded from above by
    the quantity
    $ \const \sqrt{ \gamma_{ n + 1 } } $
    where $ \const \in (0,\infty) $
    is a constant and where
    $ ( \gamma_n )_{ n \in \N } \subseteq (0,\infty) $
    is the sequence of learning
    rates;
    cf., \eg, \cite[Theorem~2.5]{DereichJentzen2024arXiv_Adam},
    \cite[Theorem~1.1]{JentzenKuckuckNeufeld2021},
    and
    \cite[Theorem~1.1]{MR4055054}.
    Further results that establish
    convergence/upper bounds for
    the strong optimization error
    for accelerated or adaptive gradient
    based optimization methods
    for deterministic or stochastic \OPs\ can,
    \eg,
    be found in
    \cite[Theorem~5.2]{Barakat_2021_cvg},
    \cite[Section~3]{MR4298987},
	\cite[Theorem~1.7]{DereichKassing2021arXiv},
    \cite[Corollary~4.11]{DereichJentzenRiekert2024arXiv_adaptive}
    and
    \cite[Section~III]{BockWeiss2019}.

	In particular, in the related work \cite{BockWeiss2019},
	convergence rates of a slightly
	modified variant of \Adam\ has also been established,
	but without explicitly known
	and, in particular, without higher order
	convergence rates,
	which is precisely the contribution
	of \cref{theo:intro:rate} above.
	In \cite{BockWeiss2019} a slightly different
	update rule of the form
	$
	  \Theta_n
	  = \Theta_{ n - 1 }
	  - \gamma
	  \frac{ \sqrt{ 1 - \beta^n } }{ 1 - \alpha^n }
	  \frac{ \m_n }{ \sqrt{ \eps + \bbM_n } }
	$
	for $ n \in \N $
	is considered.
	Unlike the update rule for \Adam\ (see,
	\eg, \cref{theo:adam_conv_coercive:eq}
	and \cref{def:Adam}),
	this function is differentiable with respect to $ \bbM $
	around $ 0 $ and, therefore, enables the use of
	standard Lyapunov tools for dynamical systems.

    For the \emph{weak optimization error
    (distance of the evaluations of the objective function)}
    for stochastic \OPs, in turn,
    one can in general
    not expect a better/lower estimate
    than a constant multiplied by the learning rate
    in the sense that
    for every $ n \in \N $ we have that the
    distance between the objective function evaluated at the optimization process at time $ n \in \N $
    and the objective function evaluated at the minimizer of the stochastic \OP\ is bounded from above by
    the quantity
    $ \const \gamma_{ n + 1 } $;
    cf., \eg,
    \cite[Section~3.3]{Chenetal2018},
    \cite[Section~3.2.2]{GodichonBaggioni2023},
    \cite[Section~4.2]{HeLiangLiuXu2023arXiv},
    and \cite[Theorem~1.1]{MR4055054}.
	Results that establish upper bounds for the weak optimization error
    for accelerated or adaptive gradient
    based optimization methods
    for deterministic \OPs\ can,
    \eg,
    also be found in
    \cite{Nesterov83,MR2142598}.

    If the objective function
    $ \cL \in C^1( \R^{ \fd }, \R ) $
    (see \cref{eq:optimization_problem} above)
    is \emph{strongly convex}
    in the sense that there exists $ c \in (0,\infty) $
    such that the modified function
    $
      \R^{ \fd } \ni \theta \mapsto \cL( \theta ) - \frac{ c }{ 2 } \| \theta \|^2 \in \R
    $
    is convex
    (cf., \eg, \cite[Section~2.3]{GarrigosGower2023arXiv_Handbook}
    and \cite[Definition~5.7.4]{JentzenBookDeepLearning2025}),
    then the objective function satisfies also
    the coercivity property that
    $
      \lim_{ r \to \infty }
      \inf_{ v \in \R^{ \fd }, \| v \| \geq r }
      \cL( v )
      = \infty
    $
    (cf., \eg, \cite[Item (ii) in Proposition~5.7.16]{JentzenBookDeepLearning2025})
    and, therefore,
    there exists a unique minimizer
    $ \lp \in \R^{ \fd } $ of the objective function
    $ \cL $
    in the sense that
    $
      \cL( \lp ) =
      \inf_{ \theta \in \R^{ \fd } } \cL( \theta )
    $
    and, hence,
	the objective function $ \cL $ satisfies
	the coercivity-type condition
	that
	for all $ \theta \in \R^{ \fd } $ it holds that
	\begin{equation}
	\label{eq:coercivity}
	  \langle
	    \theta - \lp, ( \nabla \cL )( \theta )
	  \rangle
	  \geq c \| \theta - \lp \|^2
	\end{equation}
    (cf., \eg, \cite[Proposition~5.7.16
    and Corollary~5.7.20]{JentzenBookDeepLearning2025}).
    This coercivity-type condition, in turn,
    implies that for all $ \theta \in \R^{ \fd } $
    it holds that
    $
      \cL( \theta ) - \cL( \lp )
      \geq c \| \theta - \lp \|^2
    $
    (cf., \eg,
    \cite[Lemma~5.7.22]{JentzenBookDeepLearning2025}).
    In the case of a strongly convex
    objective function $ \cL $
    we thus have that the
    squared strong optimization error
    $
%       \mathbb{E}[
        \|
          \Theta_n - \lp
        \|^2
%       ]
    $
    is, up to a constant,
    bounded by the
    weak optimization error
    $
%       \mathbb{E}[
        | \cL( \Theta_n ) - \cL( \lp ) |
%       ]
      =
%       \mathbb{E}[
        \cL( \Theta_n ) - \cL( \lp )
%       ]
    $.

    If in addition to the strong convexity
    we have that the gradient of the objective function
    $ \cL $ is Lipschitz continuous
    in the sense that there exists $ L \in \R $
    such that for all $ v, w \in \R^{ \fd } $
    it holds that
    $
      \| ( \nabla \cL )( v ) - ( \nabla \cL )( w ) \|
      \leq
      L \| v - w \|
    $ (in the context of \OP\ this
    Lipschitz property is often also referred to
    as $ L $-smoothness \cite[Definition~2.24]{GarrigosGower2023arXiv_Handbook}),
    then the above estimate also holds in the other
    direction in the sense that
    for all $ \theta \in \R^{ \fd } $
    it holds that
    $
      c \| \theta - \lp \|^2
      \leq \cL( \theta ) - \cL( \lp )
      \leq L \| \theta - \lp \|^2
    $
    (cf., \eg, \cite[Lemma~2.25]{GarrigosGower2023arXiv_Handbook} and
    \cite[Lemma 5.8.6]{JentzenBookDeepLearning2025}).
    In the strongly convex and $ L $-smooth setting,
    the squared strong optimization error
    can thus, up to a constant, be estimated by
    the weak optimization error and the other way arround.

    In this context we also note that
    the lower and the upper
    bounds for the eigenvalues of the Hessian
    $ ( \operatorname{Hess} \cL )( \lp ) $
    in
    \cref{def:convergence_order_strong}
    (corresponding to local strong convexity and local
    $ L $-smoothness, respectively)
    ensure that \cref{theo:intro:rate}
    also implies that the weak optimization
    error $ | \cL( \Theta_n ) - \cL( \lp ) | $
    for $ n \in \N $
    of \Adam\ converges with
    the squared rate
    $
      (
        \sqrt{ \cond } - 1
	  )^2
	  (
        \sqrt{ \cond } + 1
	  )^{ - 2 }
    $
    (cf.\ \cref{theo:item:adam} in \cref{theo:intro:rate})
    and, analogously, for the \RMSprop\ method
    with the squared rate
    $
      (
        \cond - 1
	  )^2
	  (
        \cond + 1
	  )^{ - 2 }
    $
    (cf.\ \cref{theo:item:rmsprop} in \cref{theo:intro:rate}).
	We also refer, \eg, to
    \cite{BarakatBianchi2020}
    for results
    for deterministic \OPs\
    in the literature
    that provide upper bounds
    for weak optimization errors.
    In particular,
    \cite[Theorem~10]{BarakatBianchi2020}
    establishes under stuitable assumptions
    convergence with rates
    for the weak optimization error
    $ | \cL( \Theta_n ) - \cL( \lp ) | $
    for $ n \in \N $
    for \Adam\ but does not provide
    information on the size of the convergence
    rate of \Adam, which is precisely the
    contribution of \cref{theo:intro:rate} above.
    In our proof of \cref{theo:intro:rate}
    we also employ the findings in \cite{BarakatBianchi2020}.
    Specifically, our proof of \cref{lem:adam_conv}
    in \cref{sec:Adam} is strongly based on applications of \cite[Theorem~2 and Proposition~14]{BarakatBianchi2020}.

    Related to the weak optimization error,
    there
    are also several works
    in the literature
    \cite{ReddiKale2019}
    that establish
    upper bounds for the so-called
    \emph{regrets of the optimization algorithm},
    that is,
    that establish
    for every $ n \in \N $ upper bounds
    for the sum
    of the empirical risk functions evaluated
    at $ \Theta_1 $, $ \Theta_2 $,
    $ \dots $, $ \Theta_n $
    subtracted by the infimal evaluation
    of the sum of the empirical risk functions;
    cf., \eg, \cite[Section~2]{ReddiKale2019}.

    There are also several works in the literature
    that provide upper bounds for moments of the
    standard norm of the gradient of the objective
    function evaluated at suitable (random) times
    of the optimization process of the
    considered optimization method; cf., \eg,
    \cite[Theorems~2 and 3]{BarakatBianchi2020},
    \cite[Theorem~3.1]{Chenetal2018},
    \cite[Section~4.1]{Defossez2022},
    \cite[Section~3.2]{HeLiangLiuXu2023arXiv},
    \cite[Theorem~1]{HongLin2024},
    \cite[Section~4]{li2023convergenceadamrelaxedassumptions},
    \cite[Section~3.1]{ZhangChen2022},
	and
	\cite[Section~3]{ZouShen2019}.
    If the objective function
    $ \cL \in C^1( \R^{ \fd }, \R ) $
    (see \cref{eq:optimization_problem} above)
    is strongly convex,
    then one gets from \cref{eq:coercivity}
    and the Cauchy-Schwarz inequality that
    for all $ \theta \in \R^d $ it holds that
    $
      \| ( \nabla \cL )( \theta ) \|
      \geq c \| \theta - \lp \|
    $.
    In the situation of strongly convex \OPs\
    we thus have that upper bounds
    for the gradient also imply upper bounds
    for the strong optimization error.

    In the situation of the \Adam\ optimizer applied
    to stochastic \OPs\ we emphasize
    that \Adam\ does typically not even converge
    to the minimizer of the stochastic \OP\ or critical points of the objective function (zeros of the gradient of the objective function)
    to which \Adam\ is applied
    but instead typically converges to different points,
    that is, to zeros of the \emph{\Adam\ vector
    field} in \cite{DereichJentzen2024arXiv_Adam};
    cf.\ \cite{ReddiKale2019}.
    We also refer, \eg, to
    \cite{ZhangChen2022,CHERIDITO2021101540,jentzen2024nonconvergenceglobalminimizersadam,Lu_2020,
	gallon2022blowphenomenagradientdescent,
	DereichGraeberJentzen2024arXiv_non_convergence}
	for further lower bounds and non-convergence results
    for \Adam\ and related \GD\ based optimization methods.
    In contrast, when applied to deterministic \OPs\ we
    do have that \Adam\ converges locally
    to local minimizers of the objective function
    due
    to, \eg, \cref{theo:item:adam}
    in \cref{theo:intro:rate} above.

    In this work
    we reveal higher order convergence
    rates for \Adam\ when applied to
    deterministic \OPs\ but
    not when applied to stochastic \OPs.
    It should, however, be pointed out that
    also in the case of stochastic \OPs\ we
    expect that one can establish
    improved error estimates involving the higher order convergence rate established here for some terms
    for the optimization error bound (that may be
    dominant in the case of large mini-batch sizes
    or early steps in the optimization process)
    by using the findings of this work in combination
    with suitable small noise analyses.

	We also refer, \eg, to the survey articles 	\cite{ruder2017overviewgradientdescentoptimization,Eetal2020survey,Sun2019}
	and monographs/handbooks  \cite{Bach2024_book_MIT_press,GarrigosGower2023arXiv_Handbook,JentzenBookDeepLearning2025}
	for further references and reviews on \Adam\ and
	related gradient based optimization methods.

	\subsection{Structure of this article}
	
	The remainder of this work is
	organized as follows.
	In \cref{sec:gradient_descent}
	we study a class of \GD\ methods 
	with abstract adaptive learning rate schedules
	covering the plain vanilla standard \GD\ and
	the \RMSprop\ methods as special cases. In particular,
	in \cref{cor:rmsprop_local} in \cref{sec:gradient_descent} we prove
	\cref{theo:item:sgd,theo:item:rmsprop} in \cref{theo:intro:rate}.
	In \cref{sec:momentum} we establish under suitable assumptions 
	convergence rates for a class of momentum \GD\ methods with abstract adaptive learning rate schedules.  
	% covering the momentum \GD\ and the \Adam\ optimizers are special cases. 
	In \cref{sec:Adam}, in turn, we first establish convergence 
	without rates of convergence for the \Adam\ optimizer and, 
	thereafter, we combine such convergence results with the general convergence 
	rate results from \cref{sec:momentum} to finally 
	reveal in \cref{cor:adam-local-c2} in 
	\cref{ssec:Adam_optimal_rate} below 
	higher order convergence rates for the \Adam\ optimizer. 
	\Cref{theo:item:adam} in \cref{theo:intro:rate} 
	follows directly from \cref{cor:adam-local-c2}.

	\subsection{Notation}
	For every $ n \in \N $,
	$ x = (x_1, \dots, x_n ) \in \R^n $
	we denote by $ \dia{x} \in \R^{ n \times n } $
	the diagonal matrix with diagonal entries
	$ x_1, x_2, \dots, x_n $
	and
	we denote by
	$ \norm{x} \in \R $
	the real number given by
	$
	\norm{x} = ( \sum_{i=1}^ n \abs{x_i}^2 )^{ 1 / 2 }
	$
	(standard norm of $ x $).
	For every
	$ n \in \N $,
	$ x = (x_1, \dots, x_n) $,
	$ y = (y_1, \dots, y_n) \in \R^n
	$
	we denote by
	$ \spro{ x, y } \in \R $
	the real number given by
	$ \spro{ x, y } = \sum_{ i = 1 }^n x_i y_i $
	(standard scalar product of $ x $ and $ y $)
	and we denote by
	$ x \cp y \in \R^n $
	the vector given by
	$
	  x \cp y =
	  \dia{x} y
% 	  =
% 	  ( x_1 y_1, \dots, x_n y_n )
% 	  ( x_i y_i )_{ i = 1 }^n
	$
	(component-wise product of $ x $ and $ y $).
	For every
	$ n \in \N $,
	$
	  ( r, x ) \in
      [
	    ( \N_0 \times \R^n )
	    \cup
	    ( [0, \infty) \times [0,\infty)^n )
	    \cup
	    ( \R \times (0,\infty)^n )
	  ]
% 	  \bigr)
	$
	we denote by
	$
	  x^{ \cp r } \in \R^n
	$
	the vector given by
	$
	  x^{ \cp r } = ( ( x_1 )^r , \dots, ( x_n )^r )
	$
	(component-wise $ r $-th power of $ x $).
	For every $ n \in \N $
	we denote by
	$ \unit{n} \in \R^{ n \times n } $
	the matrix which satisfies for all $ v \in \R^n $ that
	$
	  \unit{n} v = v
	$
	(identity matrix of size $ n $)
	and
	we denote by
    $ \vone_n \in \R^n $
    the vector given by
    $ \vone_n = ( 1, 1, \dots, 1) $
    (vector consisting of ones).
	For every $ r \in \R $,
	$ n \in \N $, $ x \in \R^n $
	we denote
	$ \B_r( x ) \subseteq \R^n $
	the set given by
	$
	  \B_r( x ) = \cu{ y \in \R^n \colon \norm{ x - y } < r }
	$
	(ball around $ x $ with radius $ r $).
	For every $ m, n \in \N $,
	$ A \in \C^{ m \times n } $
	we denote by
	$
	  \norm{A} \in \R
	$
	the real number given by
	$
	  \norm{A} = ( \sum_{i=1}^m \sum_{j=1} ^n \abs{a_{ij}} ^2 )^{ 1 / 2 }
	$
	(Frobenius norm of a matrix).
	For every $ n \in \N $, $ A \in \C^{ n \times n } $
	we denote by
	$
	  \spec{A} \subseteq \C
	$
	the set given by
	\begin{equation}
		\spec{A} = \cu{\lambda \in \C \colon \br{ \exists \, v \in \C^n \backslash \cu{0} \colon A v = \lambda v }}
	\end{equation}
	(spectrum of a square matrix)
	and we denote by $ \sr( A ) \in \R $
	the real number given by
	$
	  \sr( A ) =
	  \max(
	    \cup_{ \lambda \in \spec{A} }
	    \{ \abs{ \lambda } \}
% 	    \cu{\abs{\lambda} \colon \lambda \in \spec{A} }
	  )
	$
	(spectral radius of a square matrix).

	\section{Analysis of gradient descent (GD) with adaptive learning rates}
	\label{sec:gradient_descent}

	The goal of this section is to establish
	in \cref{cor:rmsprop_local}
	in \cref{ssec:convergence_RMSprop} below
	that the \RMSprop\ optimizer locally achieves
	the same exponential rate of convergence
	as the standard \GD\ optimizer.

\subsection{GD and RMSprop optimizers}

	In this subsection we briefly recall the definition of the standard \GD\ \cite{Cauchy1847} and
	the \RMSprop\ \cite{HintonSlides} optimizers
	in connection to
	\cref{def:convergence_order_strong} above (cf., \eg, \cite[Sections~6.5.1 and 6.5.4]{MR4688424}).

	\begin{definition}[\GD]
	\label{def:GD}
	Let $ \fd \in \N $ and let
	$ \Phi_n \colon ( \R^{ \fd } )^n \to \R^{ \fd } $,
	$ n \in \N $, be functions. Then we say that
    $ ( \Phi_n )_{ n \in \N } $ is the \GD\ optimizer
    on $ \R^{ \fd } $ (we say that
    $ ( \Phi_n )_{ n \in \N } $ is the \GD\ optimizer)
    if and only if it holds for all
    $ n \in \N $, $ g_1, g_2, \dots, g_n \in \R^{ \fd } $
    that
    \begin{equation}
	  \Phi_n( g_1, g_2, \dots, g_n )
	  =
	  g_n
	  .
    \end{equation}
	\end{definition}

    \begin{definition}[\RMSprop]
	\label{def:RMSprop}
	Let $ \fd \in \N $,
	$ \beta \in [0,1) $,
	$ \eps \in (0,\infty) $
	and let
	$
      \Phi_n = ( \Phi_n^1, \dots, \Phi_n^{ \fd } )
	  \colon
	  \allowbreak
	  ( \R^{ \fd } )^n \to \R^{ \fd }
	$,
	$ n \in \N $, be functions. Then we say that
    $ ( \Phi_n )_{ n \in \N } $ is the
    $ \beta $-$ \eps $-\RMSprop\ optimizer
    on $ \R^{ \fd } $ (we say that
    $ ( \Phi_n )_{ n \in \N } $ is
    the $ \beta $-$ \eps $-\RMSprop\ optimizer)
    if and only if it holds for all
    $ n \in \N $,
    $
      g =
      (
        ( g_{ i, j } )_{ j \in \{ 1, 2, \dots, \fd \} }
      )_{ i \in \{ 1, 2, \dots, n \} } \in ( \R^{ \fd } )^n
    $,
    $ j \in \{ 1, 2, \dots, \fd \} $
    that
    \begin{equation}
    \textstyle
	  \Phi_n^j( g )
	  =
	  \bigl[
	    \eps
	    +
	    (
% 	      ( 1 - \beta^n )^{ - 1 }
	      ( 1 - \beta )
	      \sum_{ i = 1 }^n
	      \beta^{ n - i }
	      | g_{ i, j } |^2
	    )^{ \nicefrac{ 1 }{ 2 } }
	  \bigr]^{ - 1 }
	  g_{ n, j }
	  .
    \end{equation}
	\end{definition}

	\subsection{A generalization of Gelfand's spectral radius formula}
	
	In this subsection
	we establish a generalization of the classical Gelfand formula for the spectral radius;
	see, \eg, \cite[Corollary 5.6.14]{HornJohnson2017_Book}.
	
	\newcommand{\bbK}{\mathbb{K}}
	\begin{lemma} \label{lem:equivalent-norm}
		Let $d \in \N$, $\bbK \in \cu{\R, \C }$, $A \in \bbK^{ d \times d }$.
		Then for every $\delta > 0 $ there is a norm $\norm{\cdot}_B \colon \bbK^d \to [0, \infty )$
		which satisfies for all $x \in \bbK^d$ that $\norm{A x }_B \le (\sr ( A ) + \delta ) \norm{x} _B$.
		In particular, there is a sub-multiplicative matrix norm $\norm{\cdot}_N \colon \bbK^{d \times d} \to [0, \infty )$, namely $\norm{M}_N := \sup_{x \in \bbK^d \backslash \cu{0} } \frac{\norm{M x }_B}{\norm{x} _B }$,
		such that $\norm{A}_N \le \sr ( A ) + \delta $.
	\end{lemma}
	\begin{cproof}{lem:equivalent-norm}
	The statement follows, \eg, from \cite[Lemma 5.6.10]{HornJohnson2017_Book}.
	\end{cproof}

	\begin{lemma}[Generalization of Gelfand's formula]
		\label{prop:matrix_lim}
		Let $d \in \N$,
		$M, A_1, A_2, \dots \in \C^{d \times d }$ satisfy $\lim_{n \to \infty } A_n = M$,
		and let $\delta > 0$.
		Then there exists $\const > 0$
		which satisfies for all $m, n \in \N$
		that
		\begin{equation}
			\norm*{\sprod_{i=n+1}^m A_i } \le \const ( \sr ( M ) + \delta ) ^{ m - n } .
		\end{equation}
	\end{lemma}

	\begin{cproof}{prop:matrix_lim}
		Due to \cref{lem:equivalent-norm},
		there is a sub-multiplicative matrix norm $\norm{\cdot} _B \colon \C^{d \times d } \to [0, \infty )$
		which satisfies $\norm{M}_B \le \sr ( M ) + \frac{\delta }{2}$.
		Since $A_n \to M$, there exists $N \in \N$ such that for all $n \in \N \cap [ N , \infty )$
		one has $\norm{A_n} _B \le \sr ( M ) + \delta$.
		We therefore obtain for all $m, n \in \N$ with $m>n \ge N$ that
		\begin{equation}
			\norm*{\sprod_{i=n+1}^m A_i } _B \le \sprod_{i=n+1}^m \norm{A_i} _B \le (\sr ( M ) + \delta ) ^{ m - n }.
		\end{equation}
		This demonstrates for all $m, n \in \N$ with $m > n$ that
		\begin{equation}
			\norm*{\sprod_{i=n+1}^m A_i } _B
			\le
			(\sr ( M ) + \delta ) ^{-N} \br[\big]{ \sprod_{i=1}^N \max \cu{\norm{A_i} _B , 1 } }
			( \sr ( M ) + \delta ) ^{ m - n } .
		\end{equation}
		Since $\norm{\cdot}$ and $\norm{\cdot}_B$ are equivalent, the claim follows.
	\end{cproof}
	
	\subsection{Convergence of linear iterations with perturbations}
	
	We next establish a generalization of \cite[Lemma~10]{Polyak1964}.

	\begin{lemma}
		\label{prop:iter_rate}
		Let $d \in \N$,
		$A_1, A_2, \dots, M \in \R^{ d \times d}$ satisfy $\lim_{n \to \infty} A_n = M$
		and $\sr ( M ) < 1$,
		let $h_1, h_2, \dots  \colon \R^d \to \R^d$ satisfy
		\begin{equation}
			\label{prop:iter_rate:eq_hn}
			\limsup\nolimits_{ \delta \downarrow 0 } \sup \nolimits_{n \in \N} \sup \nolimits_{y \in \B_\delta ( 0 ) \backslash \cu{0} } \frac{\norm{ h_n ( y ) } }{\norm{y } }= 0,
		\end{equation}
		let $X_0, X_1, \dots \in \R^d$ satisfy
		\begin{equation}
			\forall \, n \in \N \colon X_n = A_n X_{n-1} + h_n ( X_{n-1} ),
		\end{equation}
		and assume $\lim_{n \to \infty } X_n = 0$.
		Then for every $\delta > 0$ there is $ \const > 0$ such that
		$\forall \, n \in \N \colon \norm{X_n } \le \const ( \sr ( M ) + \delta ) ^n $.
	\end{lemma}
	
	\begin{cproof}{prop:iter_rate}
		Assume wlog that $\sr ( M ) + \delta < 1$.
		By \cref{lem:equivalent-norm} there exists a norm $\norm{\cdot}_B \colon \R^\fd \to [0, \infty )$
		which satisfies for all $y \in \R^\fd$ that $\norm{M y }_B \le (\sr ( M ) + \frac{\delta }{3}) \norm{ y } _B$.
		Since $A_n \to M$, there is $N_1 \in \N$ which satisfies for all $n \in \N \cap [N_1 , \infty )$, $y \in \R^\fd$ that
		$\norm{A_n y } _B \le (\sr ( M ) + \frac{2 \delta}{3} ) \norm{y } _B$.
		Furthermore, by the assumption on the $h_n$ there is $\eta > 0$
		such that for all $n \in \N$, $y \in \B_\eta ( 0 )$ one has $\norm{h_n ( y ) } _B \le \frac{\delta }{3} \norm{y}_B$.
		
		Now take $N_2 \in \N \cap [N_1 , \infty )$
		such that for all $n \in \N \cap [N_2, \infty )$
		it holds that $X_n \in \B_\eta ( 0 )$.
		We then obtain
		\begin{equation}
			\begin{split}
				\norm{X_{N_2+1} } _B 
				&\le \norm{A_{N_2+1} X_{N_2} } _B + \norm{h_{N_2 + 1 } ( X_{N_2} ) } _B
				\le (\sr ( M ) + \tfrac{2 \delta }{3} ) \norm{X_{N_2} }_B + \tfrac{\delta }{3} \norm{X_{N_2} }_B \\
				& = (\sr ( M ) + \delta ) \norm{X_{N_2} }_B \le \norm{X_{N_2} }_B.
			\end{split}
		\end{equation}
		By induction, it follows similarly for all $n \in \N \cap [N_2, \infty )$ that $\norm{X_{n+1} } _B \le ( \sr ( M ) + \delta ) \norm{X_n}_B$.
		In particular, there is $\const >0$ such that $\forall \, n \in \N \colon \norm{X_n} _B \le \const ( \sr ( M ) + \delta ) ^n$.
		Since $\norm{\cdot}$ and $\norm{\cdot}_B$ are equivalent, the claim follows.
	\end{cproof}

	\subsection{Convergence rates for GD with convergent learning rates}
	
	\begin{lemma}[\GD\ with convergent learning rates]
		\label{prop:gd-conv}
		Let $\fd \in \N$,
		$\kappa \in (0 , \infty )$,
		$\cK \in (\kappa , \infty )$,
		let $g \in C ( \R^\fd , \R^\fd )$ be differentiable at $0$,
		assume $g(0)=0$ and $\spec{ ( \nabla g ) ( 0 ) } \subseteq [ \kappa , \cK ]$,
		let $\Gamma_1, \Gamma_2, \dots \in (0 , \infty ) ^\fd$,
		$\Theta_0, \Theta_1, \dots \in \R^\fd$
		satisfy for all $n \in \N$ that
		\begin{equation}
			\Theta_n = \Theta_{n-1} - \Gamma_n \odot g ( \Theta_{n-1} ) ,
		\end{equation}
		and assume $\lim_{n \to \infty } \Gamma_n = \frac{2}{\kappa + \cK } \vone_\fd $
		and $\lim_{n \to \infty } \Theta_n = 0$.
		Then for every $\delta > 0$ there exists $\const > 0$ such that for all $n \in \N_0$ it holds that
		\begin{equation}
			\label{prop:gd-conv:eq-claim}
			\norm{\Theta_n } \le \const  \rbr[\big]{ \tfrac{\cK - \kappa}{\cK + \kappa } + \delta } ^n .
		\end{equation}
	\end{lemma}
	
	\begin{cproof}{prop:gd-conv}
		The fact that $g$ is differentiable at $0$ ensures that there exists a function $h \colon \R^\fd \to \R^\fd$ which satisfies $\forall \, x \in \R^\fd \colon g(x) = (\nabla g ) ( 0 ) x + h ( x )$
		and $\limsup_{x \to 0 } \frac{\norm{h(x) } }{\norm{x} } = 0$.
		In the following,
		for every $n \in \N$ let $A_n \in \R^\fd$ satisfy $A_n = \unit{\fd} - \dia{\Gamma_n } (\nabla g ) ( 0 ) $.
		We then obtain for all $n \in \N$ that
		\begin{equation}
			\Theta_n = A_n \Theta_{n-1} - \Gamma_n \odot h ( \Theta_{n-1} ) .
		\end{equation}
		Next, the assumption that
		$\lim_{n \to \infty } \Gamma_n = \frac{2}{\kappa + \cK } \vone_\fd $
		demonstrates that
		$\lim_{n \to \infty } A_n = \unit{\fd} - \frac{2}{\kappa + \cK } ( \nabla g ) ( 0 ) $.
		Moreover, \nobs that
		\begin{equation}
			\spec{ \unit{\fd} - \tfrac{2}{\kappa + \cK } ( \nabla g ) ( 0 ) } = \cu*{1 - \tfrac{2 \lambda}{\kappa + \cK } \colon \lambda \in \spec{ ( \nabla g ) ( 0 ) } } .
		\end{equation}
		The fact that $\spec{ ( \nabla g ) ( 0 ) } \subseteq [\kappa , \cK ]$
		hence demonstrates that
		$\spec{ \unit{\fd} - \tfrac{2}{\kappa + \cK } ( \nabla g ) ( 0 ) } \subseteq \br*{ - \tfrac{\cK - \kappa}{\cK + \kappa } , \tfrac{\cK - \kappa}{\cK + \kappa } } $.
		and, consequently, $\sr \rbr[\big]{ \unit{\fd} - \tfrac{2}{\kappa + \cK } ( \nabla g ) ( 0 ) } \le \tfrac{\cK - \kappa}{\cK + \kappa } $.
		Combining this with \cref{prop:iter_rate}
		and the assumption that $\lim_{n \to \infty } \Theta_n = 0$ establishes \cref{prop:gd-conv:eq-claim}.
	\end{cproof}

	\subsection{Local convergence of RMSprop}
	\label{ssec:convergence_RMSprop}

	\begin{prop}[Local convergence of \RMSprop]
		\label{lem:rmsprop_conv}
		Let $\fd \in \N$,
		$\eps, \kappa \in (0 , \infty )$,
		$\cK \in (\kappa , \infty )$,
		$\gamma \in (0 , \frac{2 \eps }{\cK } )$,
		$\beta \in [0, 1 ]$,
		let $\fl \in C^2 ( \R^\fd , \R )$
		satisfy $ (\nabla \fl ) ( 0 ) = 0$ and
		$\kappa \unit{\fd} \le ( \Hess \fl ) ( 0 ) \le \cK \unit{\fd}$,
		and for every $\theta \in \R^\fd$ let $\bbM^\theta , \Theta^\theta \colon \N_0 \to \R^\fd$
		satisfy for all $n \in \N$ that
		\begin{equation}
			\label{prop:rmsprop-local-c2:eq1}
			\begin{split}
				\Theta_0^\theta &= \theta , \qquad \bbM_0^\theta = 0 , \\
				\bbM_n^\theta &= \beta \bbM_{n-1}^\theta + ( 1 - \beta ) (\nabla \fl ( \Theta_{n-1} ^\theta ) )^{ \cp 2 }, \\
				\Theta_n ^\theta & = \Theta_{n-1} ^\theta - \gamma \rbr[\big]{ \eps \vone_\fd + (\bbM_n^\theta ) ^{ \odot 1/2} } ^{ \odot -1 } \odot ( \nabla \fl ) ( \Theta_{n-1}^\theta ) .
			\end{split}
		\end{equation}
		Then there exists $\eta > 0$ such that for all $\theta \in \B_\eta ( 0 )$
		it holds that
		\begin{equation}
			\lim\nolimits_{n \to \infty } (\norm{\Theta_n ^\theta } + \norm{\bbM_n^\theta } ) = 0.
		\end{equation}
	\end{prop}
	
	\begin{cproof}{lem:rmsprop_conv}
		For every $\theta \in \R^\fd$, $n \in \N$ let
		$\Gamma_n ^\theta =  \gamma \rbr[\big]{ \eps \vone_\fd + (\bbM_n^\theta ) ^{ \odot 1/2} } ^{ \odot -1 }$, 
		for every $\Gamma \in \R^\fd$ let
		\begin{equation}
			A_\Gamma = \unit{\fd} - \dia{\Gamma } (\Hess  \fl ) ( 0 ) ,
		\end{equation}
		and let $h_\Gamma \colon \R^{\fd } \to \R^{ \fd }$ satisfy for all $\theta \in \R^\fd$ that
		\begin{equation}
			h_\Gamma ( \theta )
			= \dia{\Gamma} \rbr*{ (\Hess ) \fl ( 0 )  \theta - \nabla \fl ( \theta ) } .
		\end{equation}
		We therefore have for all $\theta \in \R^\fd$, $n \in \N$ that
		\begin{equation}
			\Theta_n^\theta 
			= A_{\Gamma_n^\theta } 
			\Theta_{n-1}^\theta 
			+ h_{\Gamma_n^\theta } (\Theta_{n-1}^\theta ).
		\end{equation}
		The fact that $ ( \Hess  \fl ) ( 0 ) $ is symmetric
		with $\kappa \unit{\fd} \le ( \Hess \fl ) ( 0 ) \le \cK \unit{\fd}$
		ensures that $\spec{ ( \Hess  \fl ) ( 0 ) } \subseteq [ \kappa , \cK ]$.
		Combining this with the assumption that $\gamma \in (0 , \frac{2 \eps}{\cK } )$
		shows that
		\begin{equation}
			\spec{ A_{(\gamma / \eps ) \vone_\fd} } \subseteq \br[\big]{1 - \tfrac{\gamma}{\eps} \cK , 1 - \tfrac{\gamma}{\eps} \kappa }
			\subseteq ( -1 , 1 ),
		\end{equation}
		which implies that $\sr ( A_{(\gamma / \eps ) \vone_\fd}) < 1$.
		Since $A_\Gamma$ depends continuously on $\Gamma$,
		there exist $r \in (0, 1 )$, $\delta > 0$ and a norm $\norm{\cdot}_B \colon \R^{ \fd} \to [0 , \infty )$
		which satisfy for all $x \in \R^{ \fd }$, $\Gamma \in \B_\delta ( \frac{\gamma}{\eps} \vone_\fd )$
		that
		$\norm{A_\Gamma x } _B \le r \norm{x} _B$.
		Choose $\sigma_1 > 0$ such that for all
		$w \in \B_{\sigma_1} ( 0 ) \cap [0, \infty )^\fd  $
		we have that
		\begin{equation}
			\norm[\big]{\gamma \rbr[\big]{ \eps \vone_\fd + w ^{ \odot 1/2} } ^{ \odot -1 }  - \gamma \eps^{-1} \vone_\fd  } < \delta .
		\end{equation}
		Next, since $\fl \in C^2 ( \R^\fd , \R )$ and $ ( \nabla \fl ) ( 0 ) = 0$,
		there exists $\sigma _2 > 0$ such that for all $\theta \in \R^\fd$ with
		$\norm{\theta } _B < \sigma_2$ and all $\Gamma \in \B_\delta ( \frac{\gamma}{\eps} \vone_\fd )$
		it holds that $\norm{h_\Gamma (\theta ) } _B \le \frac{1-r}{2} \norm{ \theta }_B$
		and $\norm{(\nabla \fl )( \theta ) ^{ \cp 2 } } < \sigma_1 $.
		In addition, since the norms $\norm{\cdot}$, $\norm{\cdot}_B$ are equivalent,
		there is $\eta > 0$ such that for all $\theta \in \B_\eta ( 0 ) $ one has
		\begin{equation}
			\label{prop:rmsprop-local-c2:eq-ass}
			\norm{ \Theta_n^\theta } _B < \min \cu{\sigma_1, \sigma_2 }.
		\end{equation}
		Now let $\theta \in \B_\eta ( 0 )$. We show by induction that for all $n \in \N_0 $ we have
		\begin{equation}
			\label{prop:rmsprop-local-c2:eq-induct}
			\norm{\bbM_n^\theta } < \sigma_1, 
			\qquad \norm{ \Theta_n^\theta } _B < \sigma_2, 
			\qqandqq \norm{\Theta_{n}^\theta } _ B \indicator{n > 0 } \le \tfrac{1+r}{2} \norm{ \Theta_{n-1}^\theta } _B.
		\end{equation}
		In the base case $n=0$, \cref{prop:rmsprop-local-c2:eq-induct} follows from \cref{prop:rmsprop-local-c2:eq-ass}.
		Now assume that $n>0$ and that \cref{prop:rmsprop-local-c2:eq-induct} holds for $n-1$.
		\cref{prop:rmsprop-local-c2:eq1} demonstrates that
		\begin{equation}
			\begin{split}
				\norm{\bbM_n^\theta } 
				& \le \beta \norm{\bbM_{n-1}^\theta} + ( 1 - \beta ) \norm{(\nabla \fl ) ( \Theta_{n-1}^\theta ) ^{ \cp 2 } }
				\le \max \cu[\big]{\norm{\bbM_{n-1}^\theta} , \norm{(\nabla \fl ) ( \Theta_{n-1}^\theta ) ^{ \cp 2 } } } < \sigma_1 .
			\end{split}
		\end{equation}
		By definition of $\sigma_1 $, this implies that $\norm{\Gamma^\theta_n - \frac{\gamma}{\eps} \vone_\fd } < \delta$.
		Therefore, we obtain that
		\begin{equation}
			\begin{split}
				\norm{ \Theta_n^\theta } _B
				& = \norm{ A_{\Gamma^\theta_{n} }  \Theta_{n-1}^\theta + h_{\Gamma_n^\theta } ( \Theta_{n-1}^\theta ) } _B
				\le \norm{ A_{\Gamma^\theta_n }  \Theta_{n-1}^\theta }_B
				+ \norm{ h_{\Gamma_n^\theta } ( \Theta_{n-1}^\theta )  }_B \\
				& \le r \norm{ \Theta_{n-1}^\theta } _B
				+ \tfrac{1-r}{2} \norm{ \Theta_{n-1}^\theta } _B
				= \tfrac{1+r}{2} \norm{ \Theta_{n-1}^\theta } _B 
				< \sigma_2 .
			\end{split}
		\end{equation}
		Hence, \cref{prop:rmsprop-local-c2:eq-induct} also holds for $n$.
		This implies that $\lim_{n \to \infty } \Theta_n^\theta = 0$ and, consequently,
		$\lim_{n \to \infty } \bbM_n^\theta = 0$.
	\end{cproof}

	\begin{prop}[Local exponential convergence of \RMSprop]
	\label{lem:rmsprop_conv2}
	Let $\fd \in \N$,
	$\delta, \eps, \kappa \in (0 , \infty )$,
	$\cK \in (\kappa , \infty )$,
	$\gamma \in (0 , \frac{2 \eps }{\cK } )$,
	$\beta \in [0, 1 ]$,
	let $\fl \in C^2 ( \R^\fd , \R )$
	satisfy $ (\nabla \fl ) ( 0 ) = 0$ and
	$\kappa \unit{\fd} \le ( \Hess \fl ) ( 0 ) \le \cK \unit{\fd}$,
	and for every $\theta \in \R^\fd$ let $\bbM^\theta , \Theta^\theta \colon \N_0 \to \R^\fd$
	satisfy for all $n \in \N$ that
	\begin{equation}
		\label{prop:rmsprop-conv2:eq1}
		\begin{split}
			\Theta_0^\theta &= \theta , \qquad \bbM_0^\theta = 0 , \\
			\bbM_n^\theta &= \beta \bbM_{n-1}^\theta + ( 1 - \beta ) (\nabla \fl ( \Theta_{n-1} ^\theta ) )^{ \cp 2 }, \\
			\Theta_n ^\theta & = \Theta_{n-1} ^\theta - \gamma \rbr[\big]{ \eps \vone_\fd + (\bbM_n^\theta ) ^{ \odot 1/2} } ^{ \odot -1 } \odot ( \nabla \fl ) ( \Theta_{n-1}^\theta ) .
		\end{split}
	\end{equation}
	Then there exist $\eta, \const > 0$ such that for all $\theta \in \B_\eta ( 0 )$
	it holds that
	\begin{equation}
		\forall \, n \in \N \colon \norm{\Theta_n^\theta} \le \const \rbr*{ \max \cu{ \abs{ 1 - \tfrac{\gamma \kappa}{\eps} } , \abs{ 1- \tfrac{\gamma \cK}{\eps} } } + \delta } ^ n.
	\end{equation}
\end{prop}

\begin{cproof}{lem:rmsprop_conv2}
	Throughout this proof let $r = \max \cu{ \abs{ 1 - \tfrac{\gamma \kappa}{\eps} } , \abs{ 1- \tfrac{\gamma \cK}{\eps} } } \in (0 , 1 )$, assume without loss of generality that $r + \delta < 1$,
	for every $\theta \in \R^\fd$, $n \in \N$ let
	$\Gamma_n ^\theta =  \gamma \rbr[\big]{ \eps \vone_\fd + (\bbM_n^\theta ) ^{ \odot 1/2} } ^{ \odot -1 }$, 
	for every $\Gamma \in \R^\fd$ let
	\begin{equation}
		A_\Gamma = \unit{\fd} - \dia{\Gamma } (\Hess  \fl ) ( 0 ) ,
	\end{equation}
	and let $h_\Gamma \colon \R^{\fd } \to \R^{ \fd }$ satisfy for all $\theta \in \R^\fd$ that
	\begin{equation}
		h_\Gamma ( \theta )
		= \dia{\Gamma} \rbr*{ (\Hess ) \fl ( 0 )  \theta - \nabla \fl ( \theta ) } .
	\end{equation}
	We therefore have for all $\theta \in \R^\fd$, $n \in \N$ that
	\begin{equation}
		\Theta_n^\theta 
		= A_{\Gamma_n^\theta } 
		\Theta_{n-1}^\theta 
		+ h_{\Gamma_n^\theta } (\Theta_{n-1}^\theta ).
	\end{equation}
	The fact that $ ( \Hess  \fl ) ( 0 ) $ is symmetric
	with $\kappa \unit{\fd} \le ( \Hess \fl ) ( 0 ) \le \cK \unit{\fd}$
	ensures that $\spec{ ( \Hess  \fl ) ( 0 ) } \subseteq [ \kappa , \cK ]$.
	Combining this with the assumption that $\gamma \in (0 , \frac{2 \eps}{\cK } )$
	shows that
	\begin{equation}
		\spec{ A_{(\gamma / \eps ) \vone_\fd} } \subseteq \br[\big]{1 - \tfrac{\gamma}{\eps} \cK , 1 - \tfrac{\gamma}{\eps} \kappa }
		\subseteq ( -1 , 1 ),
	\end{equation}
	which implies that $\sr ( A_{(\gamma / \eps ) \vone_\fd}) \le r$.
	Since $A_\Gamma$ depends continuously on $\Gamma$,
	there exist $\sigma_1 > 0$ and a norm $\norm{\cdot}_B \colon \R^{ \fd} \to [0 , \infty )$
	which satisfy for all $x \in \R^{ \fd }$, $\Gamma \in \B_{\sigma_1} ( \frac{\gamma}{\eps} \vone_\fd )$
	that
	$\norm{A_\Gamma x}_B \le (r + \frac{\delta}{2} ) \norm{x} _B$.
	Choose $\sigma_2 > 0$ such that for all
	$w \in \B_{\sigma_2} ( 0 ) \cap [0, \infty )^\fd  $
	we have that
	\begin{equation}
		\norm[\big]{\gamma \rbr[\big]{ \eps \vone_\fd + w ^{ \odot 1/2} } ^{ \odot -1 }  - \gamma \eps^{-1} \vone_\fd  } < \sigma_1 .
	\end{equation}
	Next, since $\fl \in C^2 ( \R^\fd , \R )$ and $ ( \nabla \fl ) ( 0 ) = 0$,
	there exists $\sigma _3 > 0$ such that for all $\theta \in \R^\fd$ with
	$\norm{\theta } _B < \sigma_3$ and all $\Gamma \in \B_{\sigma_1} ( \frac{\gamma}{\eps} \vone_\fd )$
	it holds that $\norm{h_\Gamma (\theta ) } _B \le \frac{\delta}{2} \norm{ \theta }_B$
	and $\norm{(\nabla \fl )( \theta ) ^{ \cp 2 } } < \sigma_2 $.
	In addition, since the norms $\norm{\cdot}$, $\norm{\cdot}_B$ are equivalent,
	there is $\eta > 0$ such that for all $\theta \in \B_\eta ( 0 ) $ one has
	\begin{equation}
		\label{prop:rmsprop-conv2:eq-ass}
		\norm{ \Theta_n^\theta } _B < \min \cu{\sigma_2, \sigma_3 }.
	\end{equation}
	Now let $\theta \in \B_\eta ( 0 )$. We show by induction that for all $n \in \N_0 $ we have
	\begin{equation}
		\label{prop:rmsprop-conv2:eq-induct}
		\norm{\bbM_n^\theta } < \sigma_2, 
		\qquad \norm{ \Theta_n^\theta } _B < \sigma_3, 
		\qqandqq \norm{\Theta_{n}^\theta } _ B \indicator{n > 0 } \le (r + \delta ) \norm{ \Theta_{n-1}^\theta } _B.
	\end{equation}
	In the base case $n=0$, \cref{prop:rmsprop-conv2:eq-induct} follows from \cref{prop:rmsprop-conv2:eq-ass}.
	Now assume that $n>0$ and that \cref{prop:rmsprop-conv2:eq-induct} holds for $n-1$.
	\cref{prop:rmsprop-conv2:eq1} demonstrates that
	\begin{equation}
		\begin{split}
			\norm{\bbM_n^\theta } 
			& \le \beta \norm{\bbM_{n-1}^\theta} + ( 1 - \beta ) \norm{(\nabla \fl ) ( \Theta_{n-1}^\theta ) ^{ \cp 2 } }
			\le \max \cu[\big]{\norm{\bbM_{n-1}^\theta} , \norm{(\nabla \fl ) ( \Theta_{n-1}^\theta ) ^{ \cp 2 } } } < \sigma_2 .
		\end{split}
	\end{equation}
	By definition of $\sigma_2 $, this implies that $\norm{\Gamma^\theta_n - \frac{\gamma}{\eps} \vone_\fd } < \sigma_1$.
	Therefore, we obtain that
	\begin{equation}
		\begin{split}
			\norm{ \Theta_n^\theta } _B
			& = \norm{ A_{\Gamma^\theta_n }  \Theta_{n-1}^\theta + h_{\Gamma_n^\theta } ( \Theta_{n-1}^\theta ) } _B
			\le \norm{ A_{\Gamma^\theta_n }  \Theta_{n-1}^\theta }_B
			+ \norm{ h_{\Gamma_n^\theta } ( \Theta_{n-1}^\theta )  }_B \\
			& \le \rbr*{ r + \tfrac{\delta}{2} } \norm{ \Theta_{n-1}^\theta } _B
			+ \tfrac{\delta}{2} \norm{ \Theta_{n-1}^\theta } _B
			= (r + \delta) \norm{ \Theta_{n-1}^\theta } _B 
			< \sigma_3 .
		\end{split}
	\end{equation}
	Hence, \cref{prop:rmsprop-conv2:eq-induct} also holds for $n$.
	This shows for all $n \in \N$ that
	\begin{equation}
		\norm{\Theta_n^\theta } _B \le (r + \delta ) ^n \norm{\Theta_0^\theta} _B \le \sigma_2 ( r + \delta ) ^{n}.
	\end{equation}
Since $\norm{\cdot} $ and $\norm{\cdot}_B$ are equivalent, the claim follows.
\end{cproof}

	\begin{cor}[Local convergence rate for \RMSprop]
		\label{cor:rmsprop_local}
		Let $ \fd \in \N $, $ \eps, \delta, \kappa \in (0, \infty) $,
		$ \cK \in ( \kappa, \infty ) $, $ \beta \in [0, 1] $,
		$ \fl \in C^2( \R^{ \fd }, \R ) $,
		$ \lp \in \R^\fd $
		satisfy $ ( \nabla \fl )( \lp ) = 0 $ and
		$
		\kappa \unit{ \fd }
		\preceq
		%   \le
		( \Hess \fl ) ( \lp )
		%   \le
		\preceq
		\cK \unit{ \fd }
		$,
		and for every $ \theta \in \R^\fd $
		let $ \bbM^\theta \colon \N_0 \to \R^\fd $ and
		$ \Theta^\theta \colon \N_0 \to \R^\fd$
		satisfy for all $ n \in \N $ that
		\begin{equation}
% 			\label{cor:rmsprop-local-c2:eq1}
				\Theta_0^\theta = \theta ,
				\qquad
				\bbM_0^\theta = 0 ,
				% \\
				\qquad
				\bbM_n^\theta =
				\beta \bbM_{n-1}^\theta + ( 1 - \beta )
				\bigl[ (\nabla \fl )( \Theta_{n-1} ^\theta ) \bigr]^{ \cp 2 },
		\end{equation}
		\begin{equation}
				\Theta_n^\theta
				= \Theta_{n-1} ^\theta - \tfrac{2 \eps}{\kappa + \cK }
				\rbr[\big]{ \eps \vone_\fd + (\bbM_n^\theta ) ^{ \odot 1/2} } ^{ \odot ( - 1 ) } \odot ( \nabla \fl ) ( \Theta_{n-1}^\theta ) .
		\end{equation}
		Then there exist $ \eta, \const \in (0,\infty) $
		such that for all $ \theta \in \B_\eta ( \lp ) $,
		%with $\norm{ \theta - \lp } < \eta$)
         $ n \in \N $ it holds that
		\begin{equation}
			\label{eq:RMSprop_final_statement}
			\norm{ \Theta_n^\theta - \lp }
			\le
			\const
			\rbr[\big]{ \tfrac{\cK - \kappa}{\cK + \kappa } + \delta }^n .
		\end{equation}
	\end{cor}
	
	\begin{cproof}{cor:rmsprop_local}
		\Nobs that \cref{eq:RMSprop_final_statement}
		follows from \cref{lem:rmsprop_conv2} (applied with $\gamma \with \frac{2 \eps }{\kappa + \cK } $).
	\end{cproof}

	\begin{remark}
		\cref{cor:rmsprop_local}
		also includes the standard \GD\ method as a special case
		by choosing $ \eps = 1 $ and $ \beta = 0 $
		in which case the learning rates agree with the constant
		$ \frac{ 2 }{ \kappa + \cK } $.
	\end{remark}

	\section{Analysis of momentum with adaptive learning rates}
	\label{sec:momentum}

    \subsection{Momentum optimizer}

    In this subsection we briefly recall the
    definition of the momentum \GD\
    optimizer \cite{Polyak1964}
    in connection to \cref{def:convergence_order_strong} above (cf., \eg, \cite[Section~6.5.2]{MR4688424}).

    \begin{definition}[Momentum \GD]
    \label{def:momentum}
    Let $ \fd \in \N $, $ \alpha \in [0,1) $
    and let
    $
      \Phi_n = ( \Phi_n^1, \dots, \Phi_n^{ \fd } )
      \colon ( \R^{ \fd } )^n
      \allowbreak
      \to \R^{ \fd }
    $,
    $ n \in \N $, be functions. Then we say that
    $ ( \Phi_n )_{ n \in \N } $ is the
    $ \alpha $-momentum \GD\ optimizer
    on $ \R^{ \fd } $ (we say that
    $ ( \Phi_n )_{ n \in \N } $ is
    the $ \alpha $-momentum \GD\ optimizer)
    if and only if it holds for all
    $ n \in \N $,
    $
      g =
      (
        ( g_{ i, j } )_{ j \in \{ 1, 2, \dots, \fd \} }
      )_{ i \in \{ 1, 2, \dots, n \} } \in ( \R^{ \fd } )^n
    $,
    $ j \in \{ 1, 2, \dots, \fd \} $
    that
    \begin{equation}
    \textstyle
	  \Phi_n^j( g )
	  =
	  ( 1 - \alpha )
	  \bigl[
	  \sum\nolimits_{ i = 1 }^n
	  \alpha^{ n - i }
	  g_{ i, j }
	  \bigr]
	  .
    \end{equation}
	\end{definition}

	\subsection{Properties of matrix eigenvalues}
	
	\begin{lemma}
		\label{lem:block-matrix-spec}
		Let $d \in \N$,
		$\alpha \in (0 , 1 )$, $\gamma \in (0 , \infty )$,
		let $H \in \R^{d \times d}$ be symmetric,
		let $\lambda = (\lambda_1, \dots, \lambda_d ) \in \R^d$ satisfy
		$\spec{H} = \cu{\lambda_1, \ldots, \lambda_d }$,
		for every $i \in \cu{1, 2, \dots, d }$ let
		\begin{equation}
			\mu_{\pm}^{(i)} = \frac{1 + \alpha - \gamma \lambda _i (1 - \alpha)}{2} \pm 
			\sqrt{ \rbr*{\frac{1 + \alpha - \gamma \lambda_i (1 - \alpha)}{2} } ^2 - \alpha },
		\end{equation}
		let $\kappa = \min\cu{\lambda_i}$ and $\cK = \max \cu{\lambda_i }$,
		and let $A \in \R^{2d \times 2d}$ satisfy
		\begin{equation}
			A = \begin{bmatrix}
				\unit{d} - (1 - \alpha) \gamma H & - \alpha \gamma \unit{d} \\
				(1 - \alpha) H & \alpha \unit{d}.
			\end{bmatrix}
		\end{equation}
		Then
		\begin{enumerate}[label=(\roman*)]
			\item \label{lem:block-matrix-spec:item1}
			it holds that $\spec{A} = \bigcup_{i=1} ^d \cu{\mu_{\pm}^{(i) } }$,
			\item \label{lem:block-matrix-spec:item2}
			it holds that $\sr ( A ) < 1$ if and only if $\kappa > 0$ and $\gamma \cK < \frac{2 ( 1 + \alpha ) }{ 1 - \alpha }$,
			and
			\item \label{lem:block-matrix-spec:item3}
			if $\kappa > 0$, $\gamma \le \frac{1}{\sqrt{\kappa \cK } }$ and $\alpha = \rbr[\big]{\frac{1 - \gamma \kappa}{1 + \gamma \kappa } } ^2 $
			then $\sr ( A ) = \sqrt{\alpha }$. 
		\end{enumerate}
	\end{lemma}
	
	\begin{cproof}{lem:block-matrix-spec}
		The fact that $H$ is symmetric ensures that there exists an orthogonal matrix $Q \in \R^{d \times d }$ which satisfies $Q^{-1} H Q = \dia{\lambda}$.
		Hence, we obtain that
		\begin{equation}
			\begin{bmatrix}
				Q^{-1} & 0 \\
				0 & Q^{-1}
			\end{bmatrix}
			A \begin{bmatrix}
				Q & 0 \\
				0 & Q
			\end{bmatrix}
			= \begin{bmatrix}
				\unit{d} - (1 - \alpha) \gamma \dia{\lambda} & - \alpha \gamma \unit{d} \\
				(1 - \alpha ) \dia{\lambda} & \alpha \unit{d} 
			\end{bmatrix} .
		\end{equation}
		Furthermore, the fact that diagonal matrices commute and, e.g., \cite[Lemma 6.3.10]{JentzenBookDeepLearning2025}
		demonstrate for all $t \in \R$ that
		\begin{equation}
			\begin{split}
				& \det \rbr*{ \begin{bmatrix}
						\unit{d} - (1 - \alpha) \gamma \dia{\lambda} & - \alpha \gamma \unit{d} \\
						(1 - \alpha ) \dia{\lambda} & \alpha \unit{d} 
					\end{bmatrix} - t \unit{2d} } 
				= \det \rbr*{ \begin{bmatrix}
						(1 - t ) \unit{d} - (1 - \alpha) \gamma \dia{\lambda} & - \alpha \gamma \unit{d} \\
						(1 - \alpha ) \dia{\lambda} & ( \alpha - t ) \unit{d} 
				\end{bmatrix} } \\
				& = \det \rbr[\big]{ ((1 - t ) \unit{d} - (1 - \alpha) \gamma \dia{\lambda} ) ( \alpha - t) \unit{d}
					+ \alpha \gamma ( 1 - \alpha) \dia{\lambda } } \\
				&= \prod_{i=1}^d \rbr[\big]{(t - \alpha ) ( t - 1 + \gamma ( 1 - \alpha ) \lambda_i ) + \alpha \gamma ( 1 - \alpha ) \lambda_i } \\
				& = \prod_{i=1}^d \rbr[\big]{ t^2 - (1 + \alpha - (1 - \alpha ) \gamma \lambda_i ) t + \alpha } .
			\end{split}
		\end{equation}
		Therefore, we obtain that
		\begin{equation}
			\begin{split}
				\spec{A} &= \spec{ \begin{bmatrix}
						\unit{d} - (1 - \alpha) \gamma \dia{\lambda} & - \alpha \gamma \unit{d} \\
						(1 - \alpha ) \dia{\lambda} & \alpha \unit{d} 
				\end{bmatrix} } \\
				&= \cu[\Bigg]{ \frac{1 + \alpha - \gamma \lambda _i (1 - \alpha)}{2} \pm 
					\sqrt{ \rbr*{\frac{1 + \alpha - \gamma \lambda_i (1 - \alpha)}{2} } ^2 - \alpha } \colon i \in \cu{1, 2, \dots, d } } ,
			\end{split}
		\end{equation}
		which establishes \cref{lem:block-matrix-spec:item1}.
		Next, \nobs that for every 
		$i \in \cu{1, 2, \ldots, d }$ with 
		$\frac{1 - \sqrt{\alpha}}{1 + \sqrt{\alpha }} \le \gamma \lambda_i \le \frac{1 + \sqrt{\alpha}}{1 - \sqrt{\alpha }}$ it holds that
		$	\rbr*{\frac{1 + \alpha - \gamma \lambda_i (1 - \alpha)}{2} } ^2 - \alpha \le 0$.
		Hence, we obtain for every 
		$i \in \cu{1, 2, \ldots, d }$ with 
		$\frac{1 - \sqrt{\alpha}}{1 + \sqrt{\alpha }} \le \gamma \lambda_i \le \frac{1 + \sqrt{\alpha}}{1 - \sqrt{\alpha }}$
		that $\abs{\mu_-^{(i)} } = \abs{ \mu_+^{(i) } } = \sqrt{\alpha }$.
		
		Furthermore, for every $i \in \cu{1, 2, \ldots, d }$ with $\gamma \lambda_i > \frac{1 + \sqrt{\alpha}}{1 - \sqrt{\alpha }}$ one has $\abs{\mu_+^{(i)} } \le \abs{\mu_-^{ ( i ) } }$,
		and the inequality $\abs{\mu_-^{(i)} } < 1$ is equivalent to
		\begin{equation}
			\sqrt{ \rbr*{ \frac{1 + \alpha - \gamma \lambda_i ( 1 - \alpha ) }{2} } ^2 - \alpha } < \frac{3 + \alpha - \gamma \lambda_i ( 1 - \alpha ) }{2}.
		\end{equation}
		The right-hand side is positive iff $\gamma \lambda_i < \frac{3 + \alpha}{1 - \alpha }$,
		and in this case it is equivalent to
		\begin{equation}
			\begin{split}
				\rbr*{1 + \alpha - \gamma \lambda_i ( 1 - \alpha )} ^2 & < 4 \alpha + (3 + \alpha - \gamma \lambda_i ( 1 - \alpha ) ) ^2 \\
				\Longleftrightarrow
				(1 + \alpha )  ^2 - 2 ( 1+\alpha)(1 - \alpha ) \gamma \lambda_i 
				& < 4 \alpha + (3 + \alpha ) ^2 - 2 ( 3+ \alpha ) ( 1 - \alpha ) \gamma \lambda_i \\
				\Longleftrightarrow 4 ( 1 - \alpha ) \gamma \lambda_i &< 8 \alpha + 8 
				\Longleftrightarrow\gamma \lambda_i < \frac{2 (1+\alpha)}{1 - \alpha }.
			\end{split}
		\end{equation}
		Lastly, for every $i \in \cu{1, 2, \ldots, d }$ with $\gamma \lambda_i < \frac{1 - \sqrt{\alpha}}{1 + \sqrt{\alpha }} $ one has $\abs{\mu_-^{(i)} } \le \abs{\mu_+^{ ( i ) } }$,
		and the inequality $\abs{\mu_+^{( i ) } } < 1$ is equivalent to
		\begin{equation}
			\sqrt{ \rbr*{ \frac{1 + \alpha - \gamma \lambda_i ( 1 - \alpha ) }{2} } ^2 - \alpha } < \frac{1 - \alpha + \gamma \lambda_i ( 1 - \alpha ) }{2}.
		\end{equation}
		The right-hand side is positive iff $\gamma \lambda_i > -1$,
		and in this case it is equivalent to
		\begin{equation}
			\begin{split}
				\rbr*{ 1 + \alpha - \gamma \lambda_i ( 1 - \alpha ) } ^2  & < 4 \alpha + \rbr*{ 1 - \alpha + \gamma \lambda_i ( 1 - \alpha ) } ^2 \\
				\Longleftrightarrow
				(1+\alpha ) ^2 - 2 \gamma \lambda_i ( 1 + \alpha ) ( 1 - \alpha ) 
				& < 4 \alpha + ( 1 - \alpha ) ^2 + 2 \gamma \lambda_i ( 1 - \alpha ) ^2 \\
				\Longleftrightarrow
				0 & < 4 \gamma \lambda_i ( 1 - \alpha ) \Longleftrightarrow \lambda_ i > 0.
			\end{split}
		\end{equation}
		Hence, it holds for every $i \in \cu{1, 2, \ldots, d }$ that $\max \cu{\abs{\mu_-^{(i) } } , \abs{\mu_+^{ ( i ) } } } < 1$ if and only if $\kappa > 0$ and $ \gamma \cK < \frac{2 ( 1 + \alpha ) }{ 1 - \alpha } $.
		This establishes \cref{lem:block-matrix-spec:item2}.
		
		Finally, if
		$\kappa > 0$, $\gamma \le \frac{1}{\sqrt{\kappa \cK } }$ and $\alpha = \rbr[\big]{\frac{1 - \gamma \kappa}{1 + \gamma \kappa } } ^2 $
		then one has for every $i \in \cu{1, 2, \ldots, d }$
		that $\frac{1 - \sqrt{\alpha}}{1 + \sqrt{\alpha }} = \gamma \kappa \le \gamma \lambda_i$
		and
		$ \frac{1 + \sqrt{\alpha}}{1 - \sqrt{\alpha }} = \frac{1}{\gamma \kappa }
		\ge \sqrt{ \frac{\cK}{\kappa } } \ge \gamma \cK \ge \gamma \lambda_i $.
		Therefore, we obtain for every $i \in \cu{1, 2, \ldots, d }$
		that $\abs{\mu_-^{ ( i ) } } = \abs{ \mu_+ ^{ ( i ) } } = \sqrt{\alpha }$.
		This establishes \cref{lem:block-matrix-spec:item3}.
	\end{cproof}

	\subsection{Momentum GD with adaptive learning rates for quadratic functions}
	
	We next establish explicit exponential convergence rates for the momentum \GD\ method in the case of a quadratic objective function.
	The learning rates are allowed to depend on the step $n$ and the component of the iterates, as long as they converge to a fixed limit.
	This assumption will later be verified in the case of the \Adam\ optimizer.

			\begin{prop}[Convergence for sufficiently small learning rates]
				\label{prop:mom_gd2}
				Let $\fd \in \N$,
				$\lambda = (\lambda_1, \dots, \lambda_\fd) \in (0, \infty ) ^\fd$,
				$\kappa, \cK \in ( 0, \infty )$ satisfy
				$\kappa = \min \cu{\lambda_1, \dots, \lambda_\fd }$ and $\cK = \max \cu{\lambda_1, \dots, \lambda_\fd }$,
				let $\fl \colon \R^\fd \to \R$ satisfy for all $\theta = (\theta_1, \dots, \theta_\fd) \in \R^\fd$ that
				\begin{equation}
					\fl ( \theta ) = \textstyle \frac{1}{2} \sum_{i=1}^\fd \lambda_i \abs{\theta_i } ^2 ,
				\end{equation}	
				let $\gamma \in (0, \frac{1}{\sqrt{\kappa \cK } } ]$, $\alpha \in \R$ satisfy 
				$\alpha = \rbr[\big]{\frac{1 - \gamma \kappa}{1 + \gamma \kappa }}^2$,
				let
				$ \Gamma \colon \N \to (0,\infty)^{ \fd } $,
				$ \m \colon \N_0 \to \R^{ \fd } $,
				and
				$ \Theta \colon \N_0 \to \R^{ \fd } $
                satisfy for all $n \in \N$ that
				\begin{equation}
					\m_n = \alpha \m_{n-1} + (1 - \alpha) ( \nabla \fl ) ( \Theta_{n-1})
					\qqandqq \Theta_n = \Theta_{n-1} - \Gamma_n \cp \m_n,
				\end{equation}
				and assume $\lim_{n \to \infty} \Gamma_n = \gamma \vone_\fd$.
				Then for every $\delta > 0$ there exists $\const > 0$ such that for all $n \in \N_0$ it holds that
				\begin{equation}
					\label{eq:prop_mom_claim}
					\norm{\Theta_n } \le \const \rbr*{ \sqrt{\alpha} + \delta } ^n .
				\end{equation}
			\end{prop}
			
			This gives a faster convergence rate than GD, which converges with rate $(1 - \gamma \kappa )^n$ for sufficiently small $\gamma$ (see, e.g., \cite[Theorem 6.1.12]{JentzenBookDeepLearning2025}).
			
			\begin{cproof}{prop:mom_gd2}
				For every $n \in \N$ let $A_n \in \R^{ 2 \fd \times 2 \fd }$ satisfy
				\begin{equation}
					A_n = \begin{bmatrix}
						\unit{\fd} - (1 - \alpha ) \dia{ \Gamma_n \odot \lambda } & - \alpha \dia{\Gamma_n } \\
						(1 - \alpha ) \dia{\lambda} & \alpha \unit{\fd } .
					\end{bmatrix}
				\end{equation}
				The fact that for every $\theta \in \R^\fd$ it holds that $ ( \nabla \fl ) ( \theta ) = \lambda \odot \theta = \dia{\lambda} \theta$ implies for every $n \in \N$ that
				\begin{equation}
					\begin{pmatrix}
						\Theta_n \\ \m_n 
					\end{pmatrix}
					= A_n \begin{pmatrix}
						\Theta_{n-1} \\ \m_{n-1}
					\end{pmatrix} .
				\end{equation}
				This and induction imply for all $n \in \N$ that
				\begin{equation}
					\begin{pmatrix}
						\Theta_n \\ \m_n 
					\end{pmatrix}
					= A_n A_{n-1} \cdots A_1 \begin{pmatrix}
						\Theta_0 \\ \m_0 
					\end{pmatrix}.
				\end{equation}
				Moreover, the assumption that $\lim_{n \to \infty} \Gamma_n = \frac{1}{\sqrt{\kappa \cK } } \vone_\fd$ assures that
				\begin{equation}
					\lim_{n \to \infty } A_n = \begin{bmatrix}
						\unit{\fd} - \frac{1 - \alpha}{\sqrt{\kappa \cK} } \dia{\lambda} & - \frac{\alpha}{\sqrt{\kappa \cK } } \unit{\fd} \\
						(1 - \alpha ) \dia{\lambda} & \alpha \unit{\fd }
					\end{bmatrix}.
				\end{equation}
				In addition, \cref{lem:block-matrix-spec} proves that
				\begin{equation}
					\sr \rbr*{ \begin{bmatrix}
							\unit{\fd} - \frac{1 - \alpha}{\sqrt{\kappa \cK} } \dia{\lambda} & - \frac{\alpha}{\sqrt{\kappa \cK } } \unit{\fd} \\
							(1 - \alpha ) \dia{\lambda} & \alpha \unit{\fd }
					\end{bmatrix} } = \sqrt{\alpha }.
				\end{equation}
				Combining this with \cref{prop:matrix_lim}
				demonstrates that for every $\delta > 0$ there exists $\const > 0$ which satisfies for all $n \in \N$ that
				$\norm{A_n A_{n-1} \cdots A_1 } \le \const ( \sqrt{\alpha} + \delta ) ^n$.
			\end{cproof}
			
			\subsection{Convergence of momentum for locally smooth and coercive problems}
			
			We now generalize \cref{prop:mom_gd2} to the case of objective functions satisfying suitable
			local regularity properties.
			
			\begin{prop}
				\label{prop:momentum_diff_func}
				Let $\fd  \in \N$,
				$\kappa \in (0 , \infty )$, $\cK \in ( \kappa , \infty )$ ,
				let $g \in C ( \R^\fd , \R^\fd )$ be differentiable at $0$,
				assume that $g(0) = 0$, that $ ( \nabla g ) ( 0 )$ is symmetric, and $\kappa \unit{\fd} \le ( \nabla g ) ( 0 ) \le \cK \unit{\fd} $,
				let $\gamma \in (0 , \frac{1}{\sqrt{\kappa \cK } }]$, $\alpha = \rbr[\big]{ \frac{1 - \gamma \kappa}{1 + \gamma \kappa } } ^2$,
				let $\Gamma_1, \Gamma_2, \dots \in (0, \infty ) ^\fd$, $\m_0, \m_1, \dots, \Theta_0, \Theta_1, \dots \in \R^\fd$ satisfy for all $n \in \N$ that
				\begin{equation}
					\m_n = \alpha \m_{n-1} + (1 - \alpha) g ( \Theta_{n-1})
					\qqandqq \Theta_n = \Theta_{n-1} - \Gamma_n \cp \m_n,
				\end{equation}
				and assume $\lim_{n \to \infty } \Gamma_n = \gamma \vone_\fd$ and $\lim_{n \to \infty} \Theta_n = 0$.
				Then for every $\delta > 0$ there exists $\const > 0$ such that for all $n \in \N_0$ it holds that
				\begin{equation}
					\label{prop:momentum_diff_func:eq_claim}
					\norm{\Theta_n } \le \const \rbr*{ \sqrt{\alpha} + \delta } ^n .
				\end{equation}
			\end{prop}
			
			In \cref{prop:momentum_diff_func} we establish an explicit convergence rate, provided that the sequence of iterates $(\Theta_n )$ converges to the correct limit.
			Under the given assumptions on $g$ this convergence does not always hold,
			but in the case of the \Adam\ optimizer we will establish it in the next section.
			
			\begin{cproof}{prop:momentum_diff_func}
				The fact that $g$ is differentiable at $0$ ensures that there exists a function $h \colon \R^\fd \to \R^\fd$
				which satisfies $\forall \, x \in \R^\fd \colon g ( x ) = (\nabla g ) ( 0 ) x + h ( x ) $ and $\limsup_{x \to 0 } \frac{\norm{h(x)}}{\norm{x} } = 0$.
				\Nobs that for every $n \in \N$ it holds that
				$\m_n = \alpha \m_{n-1} + (1 - \alpha ) ( ( \nabla g ) ( 0 ) \Theta_{n-1} + h ( \Theta_{ n-1 } ) )$
				and
				\begin{equation}
					\Theta_n = \Theta_{n-1} - (1 - \alpha ) \Gamma_n \odot \rbr*{ ( \nabla g ) ( 0 ) \Theta_{n-1} + h ( \Theta_{ n-1 } ) }
					- \alpha \Gamma_n \odot \m_{n-1 }.
				\end{equation}
				Next, 
				for every $n \in \N$ let $A_n \in \R^{ 2 \fd \times 2 \fd }$ satisfy
				\begin{equation}
					A_n = \begin{bmatrix}
						\unit{\fd} - (1 - \alpha ) \dia{ \Gamma_n} \nabla g ( 0 ) & - \alpha \dia{\Gamma_n } \\
						(1 - \alpha ) ( \nabla g ) ( 0 ) & \alpha \unit{\fd } 
					\end{bmatrix}
				\end{equation}
				and let $H_n \colon \R^{ 2 \fd } \to \R^{ 2 \fd }$
				satisfy
				\begin{equation}
					\forall \, x, y \in \R^ \fd \colon H_n \rbr*{ \begin{pmatrix}
							x \\ y
					\end{pmatrix}}
					= \begin{pmatrix}
						- ( 1 - \alpha ) \Gamma_n \odot h ( x ) \\
						(1 - \alpha ) h ( x ) 
					\end{pmatrix} .
				\end{equation}
				We therefore obtain for all $n \in \N$ that
				\begin{equation}
					\begin{pmatrix}
						\Theta_n \\ \m_n 
					\end{pmatrix}
					= A_n 	\begin{pmatrix}
						\Theta_{n-1} \\ \m_{n-1} 
					\end{pmatrix}
					+ H _n \rbr*{ \begin{pmatrix}
							\Theta_{n-1} \\ \m_{n-1} 
					\end{pmatrix} } .
				\end{equation}
				Furthermore, the assumption that $\lim_{n \to \infty } \Gamma_n = \gamma \vone_\fd$
				demonstrates that 
				\begin{equation}
					\lim_{n \to \infty} A_n = A
					:= \begin{bmatrix}
						\unit{\fd} - (1 - \alpha ) \gamma ( \nabla g ) ( 0 ) & - \alpha \gamma \unit{\fd} \\
						(1 - \alpha ) ( \nabla g ) ( 0 ) & \alpha \unit{\fd }
					\end{bmatrix} .
				\end{equation}
				\cref{lem:block-matrix-spec} therefore demonstrates that $\sr ( A ) \le \alpha $.
				In addition, by the properties of $h$ one has
				$\limsup\nolimits_{ \delta \downarrow 0 } \sup \nolimits_{n \in \N} \sup \nolimits_{y \in \B_\delta ( 0 ) \backslash \cu{0} } \frac{\norm{ H_n ( y ) } }{\norm{y } }= 0$.
				Since $\lim_{n \to \infty } \Theta_n = 0$,
				\cref{prop:iter_rate} therefore establishes \cref{prop:momentum_diff_func:eq_claim}.
			\end{cproof}

			\section{Analysis of Adam}
			\label{sec:Adam}
			
			The goal of this section is to prove under suitable conditions that \Adam\ converges to the global minimum of the considered objective function.
			Combining this with the previous results then gives explicit exponential convergence rates.
			The two main results of this part are \cref{theo:adam_conv_coercive}, which establishes global convergence for sufficiently small learning rates,
			and \cref{cor:adam-local-c2}, which establishes local convergence with the optimal rate (for which a rather large learning rate is required).

			We first recall the definition of the \Adam\ optimizer and, thereafer,
			we recall
			in \cref{lem:adam_conv}
			that the gradients of the \Adam\ algorithm converge
			to zero for sufficiently small learning rates and objective functions
			with Lipschitz continuous gradients. 
			This has essentially been established in \cite{BarakatBianchi2020}
			and we prove \cref{lem:adam_conv} by
			applying \cite[Theorem~2 and Proposition~14]{BarakatBianchi2020}.

\subsection{Adam optimizer}

In this subsection we briefly recall the definition of the \Adam\ optimizer \cite{KingmaBa2014_Adam}
in connection
to \cref{def:convergence_order_strong} above (cf., \eg, \cite[Section~6.5.7]{MR4688424}).

    \begin{definition}[\Adam]
    \label{def:Adam}
    Let $ \fd \in \N $,
    $ \alpha, \beta \in [0,1) $,
    $ \eps \in (0,\infty) $
    and let
    $
      \Phi_n = ( \Phi_n^1, \dots, \Phi_n^{ \fd } )
      \colon
      \allowbreak
      ( \R^{ \fd } )^n
      \to \R^{ \fd }
    $,
    $ n \in \N $, be functions.
    Then we say that
    $ ( \Phi_n )_{ n \in \N } $ is the
    $ \alpha $-$ \beta $-\Adam\ optimizer
    on $ \R^{ \fd } $ (we say that
    $ ( \Phi_n )_{ n \in \N } $ is
    the $ \alpha $-$ \beta $-\Adam\ optimizer)
    if and only if it holds for all
    $ n \in \N $,
    $
      g =
      (
        ( g_{ i, j } )_{ j \in \{ 1, 2, \dots, \fd \} }
      )_{ i \in \{ 1, 2, \dots, n \} } \in ( \R^{ \fd } )^n
    $,
    $ j \in \{ 1, 2, \dots, \fd \} $
    that
    \begin{equation}
    \textstyle
	  \Phi_n^j( g )
	  =
	  \frac{ ( 1 - \alpha ) }{
	    ( 1 - \alpha^n )
	  }
% 	  ( 1 - \alpha^n )^{ - 1 }
% 	  ( 1 - \alpha )
	  \bigl[
	  \sum\nolimits_{ i = 1 }^n
	  \alpha^{ n - i }
	  g_{ i, j }
	  \bigr]
	  \bigl[
	    \eps
	    +
	    (
	      ( 1 - \beta^n )^{ - 1 }
	      ( 1 - \beta )
	      \sum_{ i = 1 }^n
	      \beta^{ n - i }
	      | g_{ i, j } |^2
	    )^{ \nicefrac{ 1 }{ 2 } }
	  \bigr]^{ - 1 }
	  .
    \end{equation}
	\end{definition}

			\subsection{Convergence of gradients to zero}
			\begin{lemma}
				\label{lem:adam_conv}
				Let $ \fd \in \N $, 
				$\cK \in (0,\infty) $,
				$\fl \in C^1 ( \R^\fd , \R )$ 
				satisfy
				\begin{equation}
					\lim\nolimits_{\norm{\theta } \to \infty } \fl ( \theta ) = \infty 
					\qqandqq \forall \, \theta , \psi \in \R^\fd \colon \norm{ (\nabla \fl ) ( \theta ) - ( \nabla \fl ) ( \psi ) } \le \cK \norm{\theta - \psi },
				\end{equation}
				let $ \eps , \gamma \in (0, \infty) $, $\alpha, \beta \in (0, 1 )$ satisfy
				$ \alpha^2 < \beta $ and
				$ \gamma < \frac{ \alpha \eps }{ \cK } $,
				and let 
				$
				\Theta,
				\m ,
				\bbM \colon \N_0 \to \R^\fd
				$ 
				satisfy for all $ n \in \N $ that
				\begin{equation}
					\label{lem:adam_conv:eq_adam}
					\begin{split}
						\textstyle 
						\m_n
						& = \alpha \m_{ n - 1 }
						+ ( 1 - \alpha ) ( \nabla \fl )( \Theta_{ n - 1 } ) , \\
						\bbM_n
						& = \beta \bbM_{ n - 1 }
						+ 
						(1 - \beta) 
						\br[\big]{( \nabla \fl ) ( \Theta_{n-1} ) } ^{ \cp 2} , \\
						\Theta_n &= 
						\Theta_{n-1}
						- \gamma ( 1 - \alpha^n ) ^{-1} 
						\rbr[\big]{ \eps \vone_\fd + ( 1 - \beta^n ) ^{ - 1/2 } ( \bbM_n ) ^{ \cp 1/2 } } ^{ \cp - 1 }
						\cp \m_n .
					\end{split}
				\end{equation}
				Then $\lim_{n \to \infty} ( \nabla \fl ) ( \Theta_n ) = 
				\lim_{n \to \infty} \bbM_n = 0$.
			\end{lemma}
			
			\begin{cproof}{lem:adam_conv}
				We will apply \cite[Theorem 2]{BarakatBianchi2020}.
				\Nobs that $\fl$ is bounded from below, continuously differentiable and $\nabla \fl$ is Lipschitz continuous with Lipschitz constant $\cK$.
				
				For every $n \in \N$ let 
				\begin{equation}
					\label{lem:adam_conv:eq_gamman}
					\Gamma_n = \gamma (1 - \alpha^n )^{-1} ( \eps \vone_\fd + (1 - \beta ^n )^{-1/2} \bbM_n^{\cp 1/2 } ) ^{ \cp - 1 } .
				\end{equation}
				We need to verify the following properties.
				\begin{enumerate}
					\item There are $N \in \N$, $\eta > \alpha$ such that $\forall \, n \in \N \cap [N, \infty ) \colon \Gamma_{n+1} \le \eta^{-1} \Gamma_n $,
					\item there are $\delta_1, \delta_2 > 0$ such that
					$C := \frac{2}{\cK} \rbr*{ 1 - \frac{(\eta - \alpha ) ^2}{2 \eta ( 1 - \alpha ) } - \frac{1 - \eta}{2(1-\alpha)} - \delta_2 } > 0$ and
					\begin{equation}
						\forall \, n \in \N \cap [N, \infty )  
						\colon \delta_1 \le \Gamma_n \le C,
					\end{equation}
				\end{enumerate}
				where the inequalities are understood componentwise.
				
				(1) Write $\bbM_n = (\bbM_n^1 , \dots , \bbM_n^\fd ) $.
				Since \cref{lem:adam_conv:eq_adam} implies for all $n \in \N$, $i \in \cu{1, \dots, \fd }$ that $\bbM_n^i \ge \beta \bbM_{n-1}^i$,
				it is enough to find some $\eta \in ( \alpha , 1 )$ such that for all $i \in \cu{1, \dots, \fd }$, $
				n \in \N \cap [N , \infty )$ one has
				\begin{equation}
					\eta ( 1 - \alpha^n ) \rbr[\Big]{ \eps + \sqrt{ ( 1 - \beta^n ) ^{-1} \bbM_n^i } }
					\le ( 1 - \alpha^{n+1} ) \rbr[\Big]{ \eps + \sqrt{ (1 - \beta^{n+1} )^{-1} \beta \bbM_n^i } }.
				\end{equation}
				Since $1 - \alpha^n \le 1 - \alpha^{n+1}$, it is sufficient if
				\begin{equation}
					\forall \, n \ge N \colon \eta \le \sqrt{ \beta ( 1 - \beta^n ) (1 - \beta^{n+1} ) ^{-1 } }
					= \sqrt{ 1 - (1 - \beta ) ( 1 - \beta^{n+1} ) ^{-1 } }.
				\end{equation}
				For $n \to \infty $ the RHS converges to $\sqrt{\beta } > \alpha$
				and thus such an $N$ exists for every $\eta \in (\alpha , \sqrt{\beta } )$.
				
				(2) 
				By \cite[Proposition 14]{BarakatBianchi2020}
				we know in particular that $(\fl ( \Theta_n ))_{n \in \N_0}$ is bounded.
				The coercivity assumption hence ensures that $\sup_{n \in \N } \norm{\Theta_n} < \infty$
				and, therefore, $\sup_{n \in \N } \norm{\bbM_n } < \infty $. 
				Combining this with \cref{lem:adam_conv:eq_gamman} implies the existence of $\delta_1$.
				For the upper bound, note that $\Gamma_n \le \frac{\gamma}{\eps ( 1 - \alpha^n ) } $.
				Furthermore,
				\begin{equation}
					\begin{split}
						1 - \frac{(\eta - \alpha ) ^2}{2 \eta ( 1 - \alpha ) } - \frac{1 - \eta}{2(1-\alpha)} 
						&= \frac{ 2 \eta ( 1 - \alpha) - \eta^2 + 2 \eta \alpha - \alpha^2 - \eta + \eta^2}{2 \eta ( 1 - \alpha ) } \\
						&= \frac{\eta - \alpha^2}{2 \eta ( 1 - \alpha ) } > \frac{\alpha - \alpha^2}{2 ( 1 - \alpha ) } = \frac{\alpha}{2} > \frac{\gamma \cK}{2 \eps} .
					\end{split}
				\end{equation}
				Hence, there exists $\delta_2 > 0$ such that $C > \frac{\gamma}{\eps}$
				and thus it holds for all sufficiently large $n$ that $\Gamma_n \le C$.
				
				We have shown that the conditions of \cite[Theorem 2]{BarakatBianchi2020} are satisfied, and hence we get that $( \nabla \fl ) (\Theta_n ) \to 0$.
				Furthermore, from the definition of $\bbM$ it is not hard to see that this implies $\bbM_n \to 0$, as desired.
			\end{cproof}

			\subsection{Convergence of Adam for quadratic problems}
			
			\begin{prop}
				\label{prop:adam_quadratic_1}
				Let $\fd \in \N$,
				$\lambda = (\lambda_1, \dots, \lambda_\fd) \in (0, \infty ) ^\fd$,
				$\kappa, \cK \in ( 0, \infty )$ satisfy
				$\kappa = \min \cu{\lambda_1, \dots, \lambda_\fd }$ and $\cK = \max \cu{\lambda_1, \dots, \lambda_\fd }$,
				let $\fl \colon \R^\fd \to \R$ satisfy for all $\theta = (\theta_1, \dots, \theta_\fd) \in \R^\fd$ that
				\begin{equation}
					\fl ( \theta ) = \textstyle \frac{1}{2} \sum_{i=1}^\fd \lambda_i \abs{\theta_i } ^2 ,
				\end{equation}
				let $\alpha, \beta \in (0, 1 )$, $\gamma, \varepsilon \in (0, \infty )$
				satisfy $\alpha = \rbr*{\frac{1 - \gamma \eps^{-1} \kappa }{ 1 + \gamma \eps^{-1} \kappa } }^2$
				and $\gamma \le \frac{\varepsilon}{\sqrt{\kappa \cK } }$,
				let $\Theta, \m, \bbM \colon \N_0 \to \R^\fd$ 
				satisfy for all $n \in \N$ that
				\begin{equation}
					\begin{split}
						\textstyle 
						\m_n
						& = \alpha \m_{ n - 1 }
						+ ( 1 - \alpha ) ( \nabla \fl )( \Theta_{ n - 1 } ) , \\
						\bbM_n
						& = \beta \bbM_{ n - 1 }
						+ 
						(1 - \beta) 
						\br[\big]{( \nabla \fl ) ( \Theta_{n-1} ) } ^{ \cp 2} , \\
						\Theta_n &= 
						\Theta_{n-1}
						- \gamma ( 1 - \alpha^n ) ^{-1} 
						\rbr[\big]{ \eps \vone_\fd + ( 1 - \beta^n ) ^{ - 1/2 } ( \bbM_n ) ^{ \cp 1/2 } } ^{ \cp - 1 }
						\cp \m_n .
					\end{split}
				\end{equation}
				and assume $ \lim_{n \to \infty} \norm{\bbM_n } = 0$.
				Then for every $\delta > 0$ there exists $\const > 0$ such that for all $n \in \N_0$ it holds that
				\begin{equation}
					\label{prop:adam_quadratic_1:eq}
					\norm{\Theta_n } \le \const \rbr*{ \sqrt{\alpha} + \delta } ^n .
				\end{equation}
			\end{prop}
			
			\begin{cproof}{prop:adam_quadratic_1}
				For every $n \in \N$ let 
				\begin{equation}
					\Gamma_n = \gamma (1 - \alpha^n )^{-1} ( \eps \vone_\fd + (1 - \beta ^n )^{-1/2} \bbM_n^{\cp 1/2 } ) ^{ \cp - 1 } .
				\end{equation}
				The assumption that $\lim_{n \to \infty} \norm{\bbM_n } = 0$ ensures that
				$\lim_{n \to \infty} \Gamma_n = \frac{\gamma}{\eps} \vone_\fd $.
				Combining this with the assumption that $\frac{\gamma}{\eps} \le \frac{1}{\sqrt{\kappa \cK } }$ and
				\cref{prop:mom_gd2} hence establishes \cref{prop:adam_quadratic_1:eq}.
			\end{cproof}

			\subsection{Convergence of Adam for smooth and coercive problems}

			\begin{cor}
				\label{cor:adam_conv_c2}
				Let $\fd \in \N$,
				$\kappa, \cK \in (0, \infty )$,
				$\fl \in C^2 ( \R^\fd, \R )$
				satisfy $ ( \nabla \fl ) ( 0 ) = 0$ and $\kappa \unit{\fd} \le ( \Hess  \fl ) ( 0 ) \le \cK \unit{\fd} $,
				let $\alpha, \beta \in (0, 1 )$, $\gamma, \varepsilon \in (0, \infty )$
				satisfy $\alpha = \rbr[\big]{\frac{1 - \gamma \eps^{-1} \kappa}{1 + \gamma \eps^{-1} \kappa }}^2$
				and $\gamma \le  \frac{\varepsilon}{\sqrt{\kappa \cK } }$,
				let $\Theta, \m, \bbM \colon \N_0 \to \R^\fd$ 
				satisfy for all $n \in \N$ that
				\begin{equation}
					\begin{split}
						\textstyle 
						\m_n
						& = \alpha \m_{ n - 1 }
						+ ( 1 - \alpha ) ( \nabla \fl )( \Theta_{ n - 1 } ) , \\
						\bbM_n
						& = \beta \bbM_{ n - 1 }
						+ 
						(1 - \beta) 
						\br[\big]{( \nabla \fl ) ( \Theta_{n-1} ) } ^{ \cp 2} , \\
						\Theta_n &= 
						\Theta_{n-1}
						- \gamma ( 1 - \alpha^n ) ^{-1} 
						\rbr[\big]{ \eps \vone_\fd + ( 1 - \beta^n ) ^{ - 1/2 } ( \bbM_n ) ^{ \cp 1/2 } } ^{ \cp - 1 }
						\cp \m_n .
					\end{split}
				\end{equation}
				and assume $\lim_{n \to \infty } \Theta_n = 0$.
				Then for every $\delta > 0$ there exists $\const > 0$ such that for all $n \in \N$ one has
				\begin{equation}
					\label{cor:adam_conv_c2:eq}
					\norm{\Theta_n } \le \const \rbr*{ \sqrt{ \alpha } + \delta } ^n .
				\end{equation}
			\end{cor}
			
			\begin{cproof}{cor:adam_conv_c2}
				For every $n \in \N$ let 
				\begin{equation}
					\Gamma_n = \gamma (1 - \alpha^n )^{-1} ( \eps \vone_\fd + (1 - \beta ^n )^{-1/2} \bbM_n^{\cp 1/2 } ) ^{ \cp - 1 } .
				\end{equation}
				The assumption that $\lim_{n \to \infty } \Theta_n = 0$ ensures that $\lim_{n \to \infty} \norm{\bbM_n } = 0$. 
				This implies that
				$\lim_{n \to \infty} \Gamma_n = \frac{\gamma}{\eps} \vone_\fd$.
				\cref{prop:momentum_diff_func} (applied with $g \with \nabla \fl$) hence establishes \cref{cor:adam_conv_c2:eq}.
			\end{cproof}

			\begin{cor}
				\label{cor:adam_conv_unique}
				Let $\fd \in \N$,
				$\kappa, \cK \in (0, \infty )$,
				let $\fl \in C^2 ( \R^\fd, \R )$ satisfy for all $\theta , \psi \in \R^\fd$ that
				\begin{equation}
					\norm{ \nabla \fl ( \theta ) - \nabla \fl ( \psi ) } \le \cK \norm{\theta - \psi } \qqandqq \spro{\theta , ( \nabla \fl ) ( \theta ) } \ge \kappa \norm{\theta } ^2 ,
				\end{equation}
				let $\alpha, \beta \in (0, 1 )$, $\gamma, \varepsilon \in (0, \infty )$
				satisfy $\alpha = \rbr[\big]{\frac{1 - \gamma \eps^{-1} \kappa}{1 + \gamma \eps^{-1} \kappa }}^2$,
				$\beta > \alpha^2$, and
				and $\gamma \le  \frac{\alpha \varepsilon}{ \cK }$,
				and let $\Theta, \m, \bbM \colon \N_0 \to \R^\fd$ 
				satisfy for all $n \in \N$ that
				\begin{equation}
					\begin{split}
						\textstyle 
						\m_n
						& = \alpha \m_{ n - 1 }
						+ ( 1 - \alpha ) ( \nabla \fl )( \Theta_{ n - 1 } ) , \\
						\bbM_n
						& = \beta \bbM_{ n - 1 }
						+ 
						(1 - \beta) 
						\br[\big]{( \nabla \fl ) ( \Theta_{n-1} ) } ^{ \cp 2} , \\
						\Theta_n &= 
						\Theta_{n-1}
						- \gamma ( 1 - \alpha^n ) ^{-1} 
						\rbr[\big]{ \eps \vone_\fd + ( 1 - \beta^n ) ^{ - 1/2 } ( \bbM_n ) ^{ \cp 1/2 } } ^{ \cp - 1 }
						\cp \m_n .
					\end{split}
				\end{equation}
				Then for every $\delta > 0$ there exists $ \const > 0$ such that for all $n \in \N$ one has
				\begin{equation}
					\norm{\Theta_n } \le \const \rbr*{ \sqrt{ \alpha } + \delta } ^n .
				\end{equation}
			\end{cor}
			
			\begin{cproof}{cor:adam_conv_unique}
				The assumptions on $\fl$ imply that $0$ is the unique global minimum of $\fl$ with
				\begin{equation}
					( \nabla \fl ) ( 0 ) = 0,
					\qquad \lim_{ \norm{\theta } \to \infty } \fl ( \theta ) = \infty , \qqandqq
					\kappa \unit{\fd} \le ( \Hess  \fl ) ( 0 ) \le \cK \unit{\fd} 
				\end{equation}
				(see, \eg, \cite[Corollary~6.1.26]{JentzenBookDeepLearning2025}).
				It therefore suffices to show that $\lim_{n \to \infty } \Theta_n = 0$,
				then the claim follows from \cref{cor:adam_conv_c2}.
				For this we observe that \cref{lem:adam_conv}
				ensures that $\lim_{n \to \infty } ( \nabla \fl ) ( \Theta_n ) = 0$.
				Combining this with the assumption that $\forall \, \theta \in \R^\fd \colon \spro{\theta , ( \nabla \fl ) ( \theta ) } \ge \kappa \norm{\theta } ^2$
				demonstrates that $\lim_{n \to \infty } \Theta_n = 0$.
			\end{cproof}

			\subsection{Faster convergence rates for sufficiently small learning rates}
			
			We next use the established convergence rates to prove that \Adam\ converges faster than \GD\ when restricted to sufficiently small learning rates.

			\begin{theorem}[\Adam\ vs.\ \GD\ with convergence rates]
			\label{theo:adam_explicit_rates}
				Let $ \fd \in \N $, 
				$
				  \lp = ( \lp_1, \dots, \lp_{ \fd } )
				$,
				$
				  \xi = ( \xi_1, \dots, \xi_{ \fd } ) \in \R^{ \fd }
				$,
				$
				  \lambda_1, \lambda_2, \dots, \lambda_{ \fd },
				  \kappa, \cK \in (0,\infty)
				$
				satisfy 
				$
				  \min_{i \in \cu{1, 2, \dots, \fd } }
				  \abs{\xi_i - \lp _i } > 0
				$
				and
				$
				  \kappa =
  				  \min\{ \lambda_1, \dots, \allowbreak \lambda_d \}
				  < \max\{ \lambda_1, \dots, \lambda_d \}
				  = \cK
				$,
				let 
				$ 
				\fl \colon \R^{ \fd } \to \R
				$ 
				satisfy for all
				$
				  \theta =
				  ( \theta_1, \dots, \theta_{ \fd } )
				  \in \R^{ \fd }
				$
				that
				\begin{equation}
				\label{eq:def_obj_function}
					\textstyle 
					\fl( \theta ) 
					= 
					\frac{1}{2} \sum_{ k = 1 }^{ \fd }
					\lambda_k ( \theta_k - \lp_k )^2 ,
				\end{equation}
				let $ \eps \in (0, 1 ) $, 
				for every $ \alpha, \beta \in (0 , 1 ) $,
				$ \gamma \in \R $ let
				$\Phi^\gamma, 
				\Theta^{ \alpha, \beta, \gamma } ,
				\m^{ \alpha, \beta, \gamma } ,
				\bbM^{ \alpha, \beta, \gamma }
				\colon \N_0 \to \R^\fd
				$ 
				satisfy for all $ n \in \N $ that
				$
				  \Phi_0^\gamma
				  = \Theta_0^{ \alpha , \beta , \gamma } = \xi
				$,
				$
				  \bbM^{ \alpha, \beta, \gamma }
				  \in [0,\infty)^{ \fd }
				$,
				and
				\begin{equation}
					\begin{split}
						\Phi_n^\gamma &= \Phi_{n-1}^\gamma - \gamma ( \nabla \fl ) ( \Phi_{n-1}^\gamma), \\
						\m_n^{ \alpha, \beta, \gamma } 
						&= \alpha \m_{ n - 1 }^{ \alpha, \beta, \gamma } 
						+ ( 1 - \alpha ) ( \nabla \fl )( \Theta_{ n - 1 }^{ \alpha, \beta, \gamma } ) 
						, \\
						\bbM_n^{ \alpha, \beta, \gamma } 
						& = \beta \bbM_{ n - 1 }^{ \alpha, \beta, \gamma } + 
						(1 - \beta) 
						\br[\big]{ (\nabla \fl ) ( \Theta_{n-1}^{ \alpha , \beta , \gamma} ) } ^{ \cp 2}, \\
						\Theta_n^{ \alpha, \beta, \gamma }
						& = 
						\Theta_{n-1}^{ \alpha, \beta, \gamma }
						- \gamma ( 1 - \alpha^n )^{-1}
						\rbr[\big]{
						  \eps \vone_\fd
						  + ( 1 - \beta^n ) ^{ - 1/2 }
						  ( \bbM_n^{ \alpha , \beta , \gamma }
						  )^{ \cp 1/2 } } ^{ \cp ( - 1 ) }
						\cp \m_n ^{ \alpha, \beta, \gamma }
						,
					\end{split}
				\end{equation}
				let $ \mlr \in (0, \frac{1}{ 4 \cK } ) $,
				$ \alpha \in \R $ satisfy
				$
				  \alpha =
				  (
				    \frac{ 1 - \mlr \kappa
				    }{ 1 + \mlr \kappa }
				  )^2
				$,
				and let
				$ \beta \in ( \alpha^2, 1 ) $,
				$ \gamma \in [0, \overline{\gamma} ] $,
				$ \delta \in (0, 1 - \mlr \kappa ) $.
				Then
				\begin{equation}
				\label{eq:statement_to_prove}
					\lim_{n \to \infty}
					\rbr*{
					  ( 1 - \mlr \kappa - \delta )^{ - n }
					  \norm{\Phi_n^\gamma - \lp }
					}
					= \infty
					\quad\text{and}\quad
					\lim_{ n \to \infty }
					\bigl(
					  \bigl( \tfrac{ 1 - \mlr \kappa }{ 1 + \mlr \kappa } + \delta \bigr)^{ \! - n }
					  \norm{ \Theta_n^{ \alpha, \beta , \eps \mlr } - \lp }
					\bigr) = 0.
				\end{equation}
			\end{theorem}
			\begin{cproof}{theo:adam_explicit_rates}
				Throughout this proof let
				$ \eta \in \R $
				satisfy
				$
				  \eta = \eps \mlr
				$.
				\Nobs that, \eg,
				\cite[Lemma~6.3.23]{JentzenBookDeepLearning2025}
				and the assumption that
				$
				  \min_{i \in \cu{1, 2, \dots, \fd } }
				  \abs{\xi_i - \lp _i } > 0
				$
				show that there exists $ c \in (0,\infty) $
				such that for all $n \in \N$ we have
				\begin{equation}
					\norm{ \Phi_n^\gamma - \lp }
					\ge
					c
					( 1 - \gamma \kappa )^n
				\end{equation}
				Hence, we obtain
				for all
				$
				  r \in ( 0, 1 - \gamma \kappa )
				$
				that
				$
				  \liminf_{ n \to \infty }
				  (
				    ( 1 - \gamma \kappa - r )^{ - n }
				    \norm{ \Phi_n^{ \gamma } - \lp }
				  )
				  = \infty
				$.
				Therefore, we obtain for all
				$
				  r \in ( \gamma \kappa, 1 )
				$
				that
				$
				  \liminf_{ n \to \infty }
				  (
				    ( 1 - r )^{ - n }
				    \norm{ \Phi_n^{ \gamma } - \lp }
				  )
				  = \infty
				$.
				The fact that
				$
				  \gamma \kappa
				  <
				  \gamma \kappa + \delta
				  \leq
				  \mlr \kappa + \delta
				  < 1
				$
				hence shows that
				\begin{equation}
				\label{eq:first_step_proof_standard_GD}
				  \lim_{ n \to \infty}
				  \bigl(
				    ( 1 - \mlr \kappa - \delta )^{ - n } \norm{\Phi_n^\gamma - \lp }
				  \bigr)
				  =
				  \lim_{ n \to \infty}
				  \bigl(
				    ( 1 - [ \mlr \kappa + \delta ] )^{ - n } \norm{\Phi_n^\gamma - \lp }
				  \bigr)
				  =
				  \infty .
				\end{equation}
				In the next step let
				$
				  \Gamma \colon \N \to \R^{ \fd }
				$
				satisfy for all $ n \in \N $ let
				\begin{equation}
					\Gamma_n  =
					\eta (1 - \alpha^n )^{-1}
					\rbr*{
					  \eps \vone_\fd
					  + (1 - \beta ^n )^{-1/2}
					  (
					    \bbM_n^{ \alpha , \beta , \eta }
					  )^{ \cp 1/2 }
					}^{ \cp ( - 1 ) } .
				\end{equation}
				\Nobs that
				the fact that
				$ 0 < \mlr \cK < \nicefrac{ 1 }{ 4 } $
				ensures that
				\begin{equation}
					\frac{ \alpha \eps }{ \cK }
					=
					\frac{ \eps }{ \cK }
					\left[
					  \frac{
					    1 - \mlr \kappa
					  }{
					    1 + \mlr \kappa
					  }
					\right]^2
					>
					\frac{ \eps }{ \cK }
					\left[
					  \frac{
					    1 - \mlr \cK
					  }{
					    1 + \mlr \cK
					  }
					\right]^2
					>
					\frac{ \eps }{ \cK }
					\left[
					  \frac{
					    (3/4)
					  }{
					    (5/4)
					  }
					\right]^2
					=
					\frac{ \eps }{ \cK }
					\left[
					  \frac{
					    3
					  }{
					    5
					  }
					\right]^2
					=
					\frac{ 9 \eps }{ 25 \cK }
% 					> \frac{ \eps }{ 4 \cK }
					> \eps \mlr
					= \eta .
				\end{equation}
				\cref{lem:adam_conv} hence shows
				that
				$
				  \lim_{n \to \infty }
				  ( \nabla \fl )(
				    \Theta_n^{ \alpha , \beta , \eta }
				  )
				  =
				  \lim_{ n \to \infty }
				  \bbM_n^{ \alpha , \beta , \eta }
				  = 0
				$.
                Combining this
                with \cref{eq:def_obj_function}
                implies that
				$
				  \lim_{n \to \infty }
				  \Theta_n^{ \alpha , \beta , \eta }
				  = \lp
				$
				and
				$
				  \lim_{ n \to \infty }
				  \bbM_n^{ \alpha , \beta , \eta }
				  = 0
				$.
				This proves that
				$
				  \lim_{n \to \infty } \Gamma_n
				  =
				  \eta
				  ( \eps \vone_{ \fd } )^{ \cp ( - 1 ) }
				  =
				  \eta \eps^{ - 1 } \vone_\fd
				  = \mlr \vone_\fd
				$.
				\cref{prop:mom_gd2}
				and the fact that
				$
				  \mlr < \frac{ 1 }{ 4 \cK }
				  <
				  \frac{ 1 }{ \cK }
				  =
				  \frac{ 1 }{ \sqrt{ \cK } \sqrt{ \cK } }
				  \leq
				  \frac{ 1 }{ \sqrt{ \kappa } \sqrt{ \cK } }
				$
				therefore
				implies
			    that
			    $
				  \lim_{ n \to \infty }
				  (
				  ( \sqrt{ \alpha } + \delta )^{ - n }
				  \norm{
				    \Theta_n^{ \alpha , \beta , \eta  } - \lp
				  }
				  ) = 0
			    $.
				Combining this with
				\cref{eq:first_step_proof_standard_GD}
				establishes \cref{eq:statement_to_prove}.
			\end{cproof}

			\begin{remark}
			We note
			in \cref{theo:adam_explicit_rates}
			that the assumption that
			$ \min_{i \in \cu{1, 2, \dots, \fd } } \abs{ \xi _i - \lp _i } > 0 $
			is only needed for the first conclusion
			$
					\lim_{n \to \infty}
					\rbr*{
					  ( 1 - \mlr \kappa - \delta )^{ - n }
					  \norm{\Phi_n^\gamma - \lp }
					}
					= \infty
			$
			in \cref{eq:statement_to_prove}
			but not for the second conclusion
			in \cref{eq:statement_to_prove}.
			Furthermore, we observe in \cref{theo:adam_explicit_rates} that
			the upper bound $ \frac{ 1 }{ 4 \cK } = \frac{ 0.25 }{ \cK } $ on $ \mlr $ could be
			replaced by $ \frac{ \sigma }{ \cK } $
			where $ \sigma \approx 0.2956 $
			is the strictly positive real solution of $ ( 1 + \sigma )^2 \sigma = ( 1 - \sigma )^2 $.
			\end{remark}

			\begin{theorem}[Smooth coercive functions]
				\label{theo:adam_conv_coercive}
				Let $\fd \in \N$, 
				$\fl \in C^2 ( \R^\fd , \R )$,
				$\lp \in \R^\fd$,
				$\kappa, \cK \in (0 , \infty )$
				satisfy for all
				$\theta , \phi \in \R^\fd$ that
				\begin{equation}
					\norm{\nabla \fl ( \theta ) - \nabla \fl ( \phi ) } \le \cK \norm{\theta - \phi }
					\qqandqq 
					\spro{\theta - \lp , \nabla \fl ( \theta ) } \ge \kappa \norm{\theta - \lp } ^2 ,
				\end{equation}
				let $\eps \in (0 , 1 )$,
				for every $ \alpha, \beta \in (0 , 1 ) $, $ \gamma \in \R $ let 
				$\Theta^{ \alpha, \beta, \gamma } ,
				\m^{ \alpha, \beta, \gamma } ,
				\bbM^{ \alpha, \beta, \gamma } \colon \N_0 \to \R^\fd
				$ 
				satisfy for all $ n \in \N $ that
				% $\m_0^{ \alpha , \beta , \gamma } = \bbM_0^{ \alpha , \beta , \gamma } = 0$,
				\begin{equation}
					\label{theo:adam_conv_coercive:eq}
					\begin{split}
						\m_n^{ \alpha, \beta, \gamma } 
						&= \alpha \m_{ n - 1 }^{ \alpha, \beta, \gamma } 
						+ ( 1 - \alpha ) ( \nabla \fl )( \Theta_{ n - 1 }^{ \alpha, \beta, \gamma } ) 
						, \\
						\bbM_n^{ \alpha, \beta, \gamma } 
						& = \beta \bbM_{ n - 1 }^{ \alpha, \beta, \gamma } + 
						(1 - \beta) 
						\br[\big]{ (\nabla \fl ) ( \Theta_{n-1}^{ \alpha , \beta , \gamma} ) } ^{ \cp 2}, \\
						\Theta_n^{ \alpha, \beta, \gamma }
						& = 
						\Theta_{n-1}^{ \alpha, \beta, \gamma }
						- \gamma ( 1 - \alpha^n )^{-1}
						\rbr[\big]{ \eps \vone_\fd + ( 1 - \beta^n ) ^{ - 1/2 } ( \bbM_n^{ \alpha , \beta , \gamma } ) ^{ \cp 1/2 } } ^{ \cp - 1 }
						\cp \m_n ^{ \alpha, \beta, \gamma }
						.
					\end{split}
				\end{equation}
				Then for all
				$\mlr \in (0, \frac{1}{4 \cK } )$
				there exist $\alpha, \beta \in (0, 1 ) $ such that for all
				% sufficiently small 
				$\delta \in (0,\infty) $
				we have
				\begin{equation}
					\inf_{\gamma \in (0, \mlr ) } \limsup_{n \to \infty} \rbr[\Big]{ \rbr[\big]{\tfrac{1 - \mlr \kappa}{1 + \mlr \kappa} + \delta} ^{-n} \norm{\Theta_n^{ \alpha, \beta , \gamma } - \lp } } = 0.
				\end{equation}
			\end{theorem}
			
			\begin{cproof}{theo:adam_conv_coercive}
				Analogous to the proof of \cref{theo:adam_explicit_rates},
				combine \cref{cor:adam_conv_unique} and \cref{lem:adam_conv}.
			\end{cproof}

			\subsection{Local convergence with optimal rate}
			\label{ssec:Adam_optimal_rate}
			
			In \cref{theo:adam_conv_coercive} we have established global exponential convergence rates for \Adam\, 
			under the restriction that the learning rate is bounded from above by $\frac{1}{4 \cK }$.
			On the other hand, to achieve the optimal rate of
			$
			  ( \sqrt{ \cK } - \sqrt{ \kappa } ) ( \sqrt{ \cK } + \sqrt{ \kappa } )^{ - 1 }
			$
			one has to use a learning rate of order $\frac{1}{\sqrt{\kappa \cK } }$,
			which can be much larger since $\frac{\kappa}{\cK }$ is usually expected to be small.
			The global convergence result in \cref{theo:adam_conv_coercive} does thus not enable us
			to conclude the optimal convergence rate for \Adam.
			
			In the following, we verify the assumption that $\lim_{n \to \infty } \Theta_n = 0$ in \cref{cor:adam_conv_c2} 
			under the assumption that the initial value $\Theta_0$ is sufficiently close to the limit point $0$.
			As a consequence, we obtain that \Adam\ converges locally with the same exponential rate as momentum \GD,
			which is generally optimal for first-order methods (cf.~Nesterov~\cite[Theorem~2.1.13]{MR2142598}).
			
			In the related work \cite{BockWeiss2019}, exponential local convergence of \Adam\ has also been established, but without an explicit rate.
			Furthermore, the authors of \cite{BockWeiss2019} consider a slightly different update rule of the form
			$\Theta_{n} = \Theta_{n-1} - \gamma \frac{\sqrt{1 - \beta^n }}{1 - \alpha^n} \frac{\m_n}{\sqrt{\eps + \bbM_n }}$.
			Unlike the update rule \cref{theo:adam_conv_coercive:eq}, this function is differentiable with respect to $\bbM$ around $0$, which enables the use of standard Lyapunov tools for dynamical systems.

			\begin{prop}
			\label{prop:adam-local-opt}
			Let $\fd \in \N$,
			$\kappa, \cK \in (0, \infty )$,
			$\fl \in C^2 ( \R^\fd, \R )$
			satisfy $ ( \nabla \fl ) ( 0 ) = 0$ and $\kappa \unit{\fd} \le ( \Hess  \fl  ) ( 0 ) \le \cK \unit{\fd} $,
			let $\alpha \in (0, 1 )$, $\beta \in [0, 1 ]$, $\gamma, \delta, \eps \in (0, \infty )$
			satisfy $\frac{\gamma}{\eps} \le \frac{1}{\sqrt{\kappa \cK } }$
			and $\alpha = \rbr{ \frac{ 1 - \gamma \eps^{-1} \kappa }{1 + \gamma \eps^{-1} \kappa } } ^2$,
			let $(a_n)_{n \in \N}, (b_n)_{n \in \N } \subseteq (0 , \infty )$
			satisfy $\lim_{n \to \infty } a_n = \lim_{n \to \infty } b_n = 1$,
			and for every $\theta \in \R^\fd$ let $\Theta^\theta , \m ^\theta , \bbM^\theta  \colon \N_0 \to \R^\fd$ 
			satisfy for all $n \in \N$ that
			\begin{equation}
				\label{prop:adam-local-opt:eq1}
				\begin{split}
					\Theta_0^\theta &= \theta , \qquad \m_0^\theta = \bbM_0^\theta = 0, \\
					\m_n^{ \theta } 
					&= \alpha \m_{ n - 1 }^{ \theta } 
					+ ( 1 - \alpha ) ( \nabla \fl )( \Theta_{ n - 1 }^{ \theta } ) 
					, \\
					\bbM_n^{ \theta } 
					& = \beta \bbM_{ n - 1 }^{ \theta } + 
					(1 - \beta) 
					\br[\big]{ (\nabla \fl ) ( \Theta_{n-1}^{ \theta } ) } ^{ \cp 2}, \\
					\Theta_n^{ \theta }
					& = 
					\Theta_{n-1}^{ \theta }
					- \gamma a_n
					\rbr[\big]{ \eps \vone_\fd + b_n ( \bbM_n^{ \theta } ) ^{ \cp 1/2 } } ^{ \cp - 1 }
					\cp \m_n ^{ \theta } .
				\end{split}
			\end{equation}
			Then there exist $\eta , \const > 0$ such that for all $\theta \in \B_\eta ( 0 )$, $n \in \N$ it holds that
			\begin{equation}
				\label{prop:adam-local-opt:eq-claim}
				\norm{\Theta_n^\theta} + \norm{\m_n^\theta } \le \const \rbr*{ \sqrt{\alpha } + \delta } ^n .
			\end{equation}
		\end{prop}
		
		\begin{cproof}{prop:adam-local-opt}
			Throughout this proof assume without loss of generality that $\sqrt{\alpha } + \delta < 1$,
			for every $\theta \in \R^\fd$, $n \in \N$ define
			\begin{equation}
				\Gamma_n ^\theta =  \gamma a_n \rbr[\big]{ \eps \vone_\fd  + b_n (\bbM_n^\theta  ) ^{\cp 1/2 } } ^{ \cp - 1 } ,
			\end{equation}
			for every $\Gamma \in \R^\fd$ let
			\begin{equation}
				A_\Gamma = \begin{bmatrix}
					\unit{\fd} - ( 1 - \alpha ) \dia{\Gamma} ( \Hess  \fl ) ( 0 ) & - \alpha \dia{\Gamma } \\
					(1 - \alpha ) (\Hess  \fl ) ( 0 ) & \alpha \unit{\fd}
				\end{bmatrix}
				\in \R^{ 2 \fd \times 2 \fd } 
			\end{equation}
			and let $h_\Gamma \colon \R^{2 \fd } \to \R^{ 2 \fd }$ satisfy for all $\theta , m \in \R^\fd$ that
			\begin{equation}
				h_\Gamma \rbr*{ \begin{pmatrix}
						\theta \\ m
				\end{pmatrix}}
				= \begin{pmatrix}
					(1 - \alpha ) \dia{\Gamma} ( (\Hess  \fl ) ( 0 )  \theta - ( \nabla \fl ) ( \theta ) ) \\
					(1 - \alpha ) ( ( \nabla \fl ) ( \theta ) - ( \Hess  \fl ) ( 0 ) \theta )
				\end{pmatrix}
			\end{equation}
			We therefore have for all $\theta \in \R^\fd$, $n \in \N$ that
			\begin{equation}
				\begin{pmatrix}
					\Theta_n^\theta \\ \m_n^\theta 
				\end{pmatrix}
				= A_{\Gamma_n^\theta } \begin{pmatrix}
					\Theta_{n-1}^\theta \\ \m_{n-1}^\theta 
				\end{pmatrix}
				+ h_{\Gamma_n^\theta } \rbr*{ \begin{pmatrix}
						\Theta_{n-1}^\theta \\ \m_{n-1}^\theta 
				\end{pmatrix} }
			\end{equation}
			The fact that $( \Hess  \fl ) ( 0 ) $ is symmetric ensures that $\spec{( \Hess  \fl ) ( 0 ) } \subseteq [ \kappa , \cK ]$.
			Hence, \cref{lem:block-matrix-spec} implies that $\sr ( A_{(\gamma / \eps ) \vone_\fd}) = \sqrt{\alpha }$.
			Since $A_\Gamma$ depends continuously on $\Gamma$,
			there exist $\nu > 0$ and a norm $\norm{\cdot}_B \colon \R^{ 2 \fd} \to [0 , \infty )$
			which satisfy for all $x \in \R^{ 2 \fd }$, $\Gamma \in \B_\nu ( \frac{\gamma}{\eps} \vone_\fd )$
			that
			$\norm{A_\Gamma x } _B \le ( \sqrt{\alpha } + \frac{\delta}{2} ) \norm{x} _B$.
			The fact that $\lim_{n \to \infty } a_n = \lim_{n \to \infty } b_n = 1$ ensures that there exist $N \in \N$, $\sigma_1 > 0$ which satisfy for all $n \in \N \cap [N, \infty )$,
			$w \in \B_{\sigma_1} ( 0 ) \cap [0, \infty )^\fd  $
			that
			\begin{equation}
				\norm*{\gamma a_n ( \eps \vone_\fd  + b_n w ^{\cp 1/2 } ) ^{ \cp - 1 }  - \gamma \eps^{-1} \vone_\fd  } < \nu .
			\end{equation}
			Furthermore, since $\fl \in C^2 ( \R^\fd , \R )$ and $( \nabla \fl ) ( 0 ) = 0$,
			there exists $\sigma _2 > 0$ such that for all $\theta , m \in \R^\fd$ with
			$\norm{(\theta , m ) } _B < \sigma_2$ and all $\Gamma \in \B_\nu ( \frac{\gamma}{\eps} \vone_\fd )$
			it holds that $\norm{h_\Gamma (\theta , m ) } _B \le \frac{\delta }{2} \norm{ (\theta , m ) }_B$
			and $\norm{(\nabla \fl ( \theta ) ) ^{ \cp 2 } } < \sigma_1 $.
			Now since the sequences $\Theta^\theta, \m^\theta, \bbM^\theta$ depend continuously on $\theta$,
			there is $\eta > 0$ such that for all $\theta \in \B_\eta ( 0 ) $, $n \in \cu{0, 1, \ldots, N }$ one has
			\begin{equation}
				\label{prop:adam-local-opt:eq-ass}
				\norm{\bbM_n^\theta } + \norm{ ( \Theta_n^\theta , \m_n^\theta ) } _B < \min \cu{\sigma_1, \sigma_2 }.
			\end{equation}
			We now show by induction that for all $n \in \N \cap [N, \infty )$, $\theta \in \B_\eta ( 0 )$ we have that
			\begin{equation}
				\label{prop:adam-local-opt:eq-induct}
				\norm{\bbM_n^\theta } < \sigma_1, 
				\qquad \norm{( \Theta_n^\theta , \m_n^\theta ) } _B < \sigma_2, 
				\qqandqq \norm{(\Theta_{n}^\theta , \m_{n}^\theta ) } _ B \indicator{n > N } \le (\sqrt{ \alpha } + \delta ) \norm{ ( \Theta_{n-1}^\theta , \m_{n-1}^\theta ) } _B.
			\end{equation}
			In the base case $n=N$, \cref{prop:adam-local-opt:eq-induct} follows from \cref{prop:adam-local-opt:eq-ass}.
			Now assume that $n>N$ and that \cref{prop:adam-local-opt:eq-induct} holds for $n-1$.
			\cref{prop:adam-local-opt:eq1} demonstrates for all $\theta \in \B_\eta ( 0 ) $ that
			\begin{equation}
				\begin{split}
					\norm{\bbM_n^\theta } 
					& \le \beta \norm{\bbM_{n-1}^\theta} + ( 1 - \beta ) \norm{(\nabla \fl ) ( \Theta_{n-1}^\theta ) ^{ \cp 2 } }
					\le \max \cu[\big]{\norm{\bbM_{n-1}^\theta} , \norm{(\nabla \fl ) ( \Theta_{n-1}^\theta ) ^{ \cp 2 } } } < \sigma_1 .
				\end{split}
			\end{equation}
			By definition of $\sigma_1 $, this implies that $\norm{\Gamma^\theta_n - \frac{\gamma}{\eps} \vone_\fd } < \nu $.
			Therefore, we obtain for all $\theta \in \B_\eta ( 0 ) $ that
			\begin{equation}
				\begin{split}
					 \norm{( \Theta_n^\theta , \m_n^\theta ) } _B
					& = \norm{ A_{\Gamma^\theta_n } ( \Theta_{n-1}^\theta , \m_{n-1}^\theta ) + h_{\Gamma_n^\theta } ( \Theta_{n-1}^\theta , \m_{n-1}^\theta ) } _B \\
					& \le \norm{ A_{\Gamma^\theta_n } ( \Theta_{n-1}^\theta , \m_{n-1}^\theta ) }_B
					+ \norm{ h_{\Gamma_n^\theta } ( \Theta_n^\theta , \m_n^\theta )  }_B \\
					& \le \rbr*{ \sqrt{ \alpha } + \tfrac{\delta}{2} } \norm{( \Theta_{n-1}^\theta , \m_{n-1}^\theta ) } _B
					+ \tfrac{\delta }{2} \norm{( \Theta_{n-1}^\theta , \m_{n-1}^\theta ) } _B \\
					& = \rbr*{ \sqrt{ \alpha } + \delta } \norm{( \Theta_{n-1}^\theta , \m_{n-1}^\theta ) } _B 
					< \sigma_2 .
				\end{split}
			\end{equation}
			Hence, \cref{prop:adam-local-opt:eq-induct} also holds for $n$.
			Consequently, we obtain for all $\theta \in \B_\eta ( 0 ) $, $n \in \N \cap [ N , \infty )$
			that
			$\norm{(\Theta_n^\theta , \m_n^\theta ) } _B \le ( \sqrt{\alpha} + \delta ) ^{n - N} \sigma_2 $.
			This and \cref{prop:adam-local-opt:eq-ass}
			demonstrate for all
			$\theta \in \B_\eta ( 0 ) $, $n \in \N $ that
			\begin{equation}
				\begin{split}
					\norm{(\Theta_n^\theta , \m_n^\theta ) } _B
					\le \br[\big]{  \sigma_2 \rbr*{ \sqrt{\alpha} + \delta } ^{-N} } \rbr*{ \sqrt{\alpha } + \delta } ^n .
				\end{split}
			\end{equation}
		The fact that $\norm{\cdot}$ and $\norm{\cdot}_B$ are equivalent therefore establishes \cref{prop:adam-local-opt:eq-claim}.
		\end{cproof}

					\begin{cor}
			\label{cor:momentum-local-c2}
			Let $\fd \in \N$, $ \lp \in \R^{ \fd } $,
			$\kappa, \cK \in (0, \infty )$,
			$\fl \in C^2 ( \R^\fd, \R )$
			satisfy $( \nabla \fl ) ( \lp ) = 0$ and $\kappa \unit{\fd} \le ( \Hess  \fl ) ( \lp ) \le \cK \unit{\fd} $,
			let $\alpha \in (0, 1 )$, $\gamma, \delta \in (0, \infty )$
			satisfy $\alpha = \rbr[\big]{\frac{1 - \gamma \kappa}{1 + \gamma \kappa }}^2$
			and $\gamma \le  \frac{1}{\sqrt{\kappa \cK } }$,
			and for every $\theta \in \R^\fd $
			let
			$ \m^\theta \colon \N_0 \to \R^\fd$
			and
			$ \Theta^\theta \colon \N_0 \to \R^{ \fd } $
			satisfy for all $n \in \N$ that
			\begin{equation}
					\m_0^\theta = 0,
					\qquad
					\m_n^{ \theta } 
					= \alpha \m_{ n - 1 }^{ \theta }
					+ ( 1 - \alpha )
					( \nabla \fl )(
					  \Theta_{ n - 1 }^{ \theta }
					)
					,
			\end{equation}
			\begin{equation}
					\Theta_0^\theta = \theta ,
					\qquad
					\text{and}
					\qquad
					\Theta_n^{ \theta }
					=
					\Theta_{n-1}^{ \theta }
					- \gamma \m_n ^{ \theta } .
			\end{equation}
			Then there exist
			$ \eta, \const \in (0,\infty) $
			such that for all
			$\theta \in \B_\eta ( \lp )$, $n \in \N$ one has
			\begin{equation}
				\norm{ \Theta_n^\theta - \lp }
				\le \const \rbr*{ \sqrt{ \alpha } + \delta } ^n .
			\end{equation}
		\end{cor}
	
	\begin{cproof}{cor:momentum-local-c2}
		This is a direct consequence of \cref{prop:adam-local-opt}
		(applied with $\eps \with 1$, $\beta \with 1$, $a_n \with 1$, $b_n \with 1$
		in the notation of \cref{prop:adam-local-opt}).
	\end{cproof}

			\begin{cor}
				\label{cor:adam-local-c2}
				Let $\fd \in \N$,
				$\kappa, \cK \in (0, \infty )$,
				$\fl \in C^2 ( \R^\fd, \R )$
				satisfy $( \nabla \fl ) ( 0 ) = 0$ and $\kappa \unit{\fd} \le ( \Hess  \fl ) ( 0 ) \le \cK \unit{\fd} $,
				let $\alpha, \beta \in (0, 1 )$, $\gamma, \varepsilon \in (0, \infty )$
				satisfy $\alpha = \rbr[\big]{\frac{1 - \gamma \eps^{-1} \kappa}{1 + \gamma \eps^{-1} \kappa }}^2$
				and $\gamma \le  \frac{\varepsilon}{\sqrt{\kappa \cK } }$,
				and for every $\theta \in \R^\fd $
				let $\Theta^\theta , \m^\theta, \bbM^\theta \colon \N_0 \to \R^\fd$ 
				satisfy for all $n \in \N$ that
				\begin{equation}
					\begin{split}
						\Theta_0^\theta &= \theta , \qquad \m_0^\theta = \bbM_0^\theta = 0, \\
						\m_n^{ \theta } 
						&= \alpha \m_{ n - 1 }^{ \theta } 
						+ ( 1 - \alpha ) ( \nabla \fl )( \Theta_{ n - 1 }^{ \theta } ) 
						, \\
						\bbM_n^{ \theta } 
						& = \beta \bbM_{ n - 1 }^{ \theta } + 
						(1 - \beta) 
						\br[\big]{ (\nabla \fl ) ( \Theta_{n-1}^{ \theta } ) } ^{ \cp 2}, \\
						\Theta_n^{ \theta }
						& = 
						\Theta_{n-1}^{ \theta }
						- \gamma ( 1 - \alpha^n )^{-1}
						\rbr[\big]{ \eps \vone_\fd + ( 1 - \beta^n ) ^{ - 1/2 } ( \bbM_n^{ \theta } ) ^{ \cp 1/2 } } ^{ \cp - 1 }
						\cp \m_n ^{ \theta } .
					\end{split}
				\end{equation}
				Then there exist $\eta , \const > 0$  such that for all
				$\theta \in \B_\eta ( 0 )$, $n \in \N$ one has
				\begin{equation}
					\norm{\Theta_n ^\theta } \le \const \rbr*{ \sqrt{ \alpha } + \delta } ^n .
				\end{equation}
			\end{cor}
			
			\begin{cproof}{cor:adam-local-c2}
				This is a direct consequence of \cref{prop:adam-local-opt}
				(applied with $a_n \with (1 - \alpha^n)^{-1}$, $b_n \with (1 - \beta^n)^{-1/2 }$).
			\end{cproof}

			\subsection*{Acknowledgements}
			
			This work has been supported by the
			Ministry of Culture
			and Science NRW as part of the Lamarr Fellow Network.
			In addition, this work has been partially funded
			by the Deutsche Forschungsgemeinschaft
			(DFG, German Research Foundation) under
			Germany's Excellence Strategy EXC 2044-390685587,
			Mathematics Münster: Dynamics-Geometry-Structure.
			Robin Graeber is gratefully acknowledged for
			several useful discussions.
			
			\bibliographystyle{acm}
			\bibliography{bibfile}
\end{document}